\documentclass[11pt]{article}

\makeatletter
\let\@fnsymbol\@arabic
\makeatother

\usepackage{appendix}
\usepackage{graphicx}
\usepackage{makecell}
\usepackage{amsmath}
\usepackage{amsfonts}
\usepackage{amssymb}
\usepackage{epsfig}

\usepackage{amsmath}
\usepackage{amsfonts}
\usepackage{amssymb}
\usepackage{dsfont}
\usepackage{mathrsfs}
\usepackage{bbm}
 \usepackage{hyperref}
\usepackage{amsmath}
\usepackage{amsfonts}
\usepackage{amssymb}
\usepackage{multirow}
\usepackage{mathrsfs}
\usepackage[dvipsnames,usenames]{color}

\usepackage{subcaption}
\usepackage{caption}
\usepackage{svg}

\usepackage{natbib}
\setcitestyle{numbers,square}
\renewcommand{\Box}{\framebox{\rule{0.3em}{0.0em}}}

\newtheorem{theorem}{Theorem}[section]
\newtheorem{theorem*}{Theorem}[subsubsection]
\newtheorem{lemma}[theorem]{Lemma}
\newtheorem{proposition}[theorem]{Proposition}
\newtheorem{example}{Example}[section]

\newtheorem{definition}[theorem]{Definition}

\newtheorem{assumption}{Assumption}
\newtheorem{corollary}[theorem]{Corollary}

\newcommand{\setd}{{ d \kern -.15em l}}
\newcommand{\hatsetd}{ d \hat{\kern -.15em l }}
\newcommand{\dd}{\mathsf {d\kern -0.07em l}} %	nested

\newcommand{\z}{{ z}}

\newcommand{\C}{{\cal C}}

\newcommand{\N}{{\cal N}}

\newcommand{\bgeqn}{\begin{eqnarray}}
\newcommand{\edeqn}{\end{eqnarray}}
\newcommand{\bgeq}{\begin{eqnarray*}}
\newcommand{\edeq}{\end{eqnarray*}}
\newcommand{\bec}{\begin{center}}
\newcommand{\enc}{\end{center}}
\newcommand{\R}{{\rm I\!R}}

\newcommand{\inmat}[1]{\mbox{\rm {#1}}}

\newcommand{\Z}{{\cal Z}}
\newcommand{\F}{{\cal F}}

\newcommand{\A}{{\cal A}}
\newcommand{\B}{{\cal B}}
\newcommand{\G}{{\cal G}}

\newcommand{\X}{{\cal X}}

\newcommand{\W}{\mathcal{W}}

\newcommand{\be}{\begin{equation}}
\newcommand{\ee}{\end{equation}}

\def\w{\omega}

\def\bbr{{\Bbb{R}}} %real numbers
\def\bbe{{\Bbb{E}}} %expectation

\renewcommand{\Box}{\hfill \rule{2.3mm}{2.3mm}}

\numberwithin{equation}{section}

\title{
Necessary Optimality Conditions for Integrated Learning and Optimization 
Problem in Contextual Optimization 
}

 \author{
 Yuan Tao\thanks{Department of Systems Engineering and Engineering Management, The Chinese University of Hong Kong. Email: ytao@se.cuhk.edu.hk}
 \,\,and\,\,
 Huifu Xu\thanks{Department of Systems Engineering and Engineering Management, The Chinese University of Hong Kong. Email: hfxu@se.cuhk.edu.hk.}
 }

\begin{document}

\maketitle

 \begin{abstract}

Integrated learning and optimization (ILO)
 is a framework in contextual optimization 
 which aims to train a predictive model for the probability distribution of the underlying 
 problem data uncertainty, 
with the goal of enhancing the quality of downstream decisions. This framework represents a new class of stochastic bilevel programs, which are extensively utilized in the literature of operations research and management science, yet remain underexplored from the perspective of optimization theory.
In this paper, we fill out the gap. 
Specifically, 
we derive the first-order necessary optimality conditions 
 in terms of Mordukhovich limiting subdifferentials.
To this end, we formulate bilevel program as a two-stage 
stochastic program  with variational inequality constraints 
when the lower-level decision-making problem is convex, and establish an optimality condition via sensitivity analysis of the second-stage value function. In the case where the lower level optimization problem is nonconvex, 
we adopt the value function approach in the literature of 
bilevel program and derive the first-order necessary conditions under stochastic partial calmness conditions. 
The derived optimality conditions are applied to several existing ILO problems in the literature.
These conditions may be used for the design of gradient-based algorithms for solving ILO problems.
 \end{abstract}

 \textbf{Keywords.} Integrated learning and optimization, Mordukhovich limiting subdifferential, stochastic bilevel program, partial calmness conditions

% {Keywords}

\section{Introduction} \label{sec:model}

In this paper, we study the first-order optimality conditions of the following integrated learning and optimization (ILO) problem
\begin{subequations}
\label{eq:ILO-1}
    \begin{eqnarray}
     % (\inmat{ILO--1})\qquad
        \min_{\theta\in\Theta, z^*(\cdot)}
        %\in \mathscr{Z}}
        && 
        % \mathbb{E}_{P_X}\left[\mathbb{E}_{P_{Y|X}}[L(z^*(\theta,x),X,Y)]\right] = 
        \mathbb{E}_{P} \left[
        %\min_{z^*(\theta,x)}
        L(z^*(X),X,Y,\theta)\right] \label{eqn:upper-level-problem_A}\\
        \inmat{s.t.} && z^*(x) \in  \arg\min_{z\in \mathcal{Z}} \mathfrak{C}(z,\theta,x) :=\mathbb{E}_{f_\theta(x)}[c(z,Y)],\; \inmat{for a.e.} \;x\in \mathcal{X}, \quad
        \label{eqn:lower-level-problem-A}
        % x \in \mathcal{X}.
\end{eqnarray}
\end{subequations}
%This is a bilevel program. At the upper level,
where 
$L:\mathbb{R}^{d_z}\times \mathbb{R}^{d_x}\times \mathbb{R}^{d_y}\times \mathbb{R}^{d_\theta} 
\rightarrow \mathbb{R}$ is a continuous function, 
$X$ and $Y$
are random vectors mapping from $(\Omega,\mathcal{F},\mathbb{P})$ to $\mathbb{R}^{d_x}$ and $ \mathbb{R}^{d_y}$ with support sets ${\cal X}$ and ${\cal Y}$, $P$ is the joint probability distribution of $(X,Y)$; 
$c: \mathbb{R}^{d_z}\times \mathbb{R}^{d_y} \rightarrow \mathbb{R}$ is a continuous function, $f_\theta(x)$ represents the predicted  probability distribution of $Y$ parameterized by $\theta\in \Theta \subset \mathbb{R}^{d_\theta}$ based on the observation of $X=x$, and
$z^*(\cdot)$ is a measurable selection from optimal solution mapping from
$\mathcal{X}$ to $\mathcal{Z}$.

This is a two-stage stochastic bilevel program: at the lower level (the second stage), given contextual information \( x \), an optimal decision \( z^*(x) \) is selected to minimize the expected cost \( c \) under the predicted distribution \( f_\theta(x) \) of \( Y \), which is parameterized by \( \theta \);
 at the upper level (the first stage), the quality of the prediction is assessed via the expected value of a loss function \( L \), which depends on the lower-level decision \( z^*(x) \) and the true distribution of \( (X, Y) \). The expectation is taken with respect to the true joint distribution \( P \). An optimal parameter \( \theta^* \) is then selected to yield the best predictive model of the conditional distribution of \( Y \) given \( X \). Problem \eqref{eq:ILO-1} originates from the following conditional stochastic optimization problem:
\begin{equation*}
    \min_{z \in \mathcal{Z}} \, \mathbb{E}_{P}[c(z, Y) \mid X = x],
\end{equation*}
in which the uncertainty in the problem data, represented by the random variable \( Y \), is statistically dependent on the contextual variable \( X \), which is observed prior to the realization of \( Y \). In practice, the relationship between \( X \) and \( Y \) is typically inferred from empirical data.
However, the conditional distribution $ P(Y \mid X = x) $ is rarely available in closed form.
A common modeling approach in contextual optimization is to approximate the underlying conditional distribution with a parametric family $ \{ f_\theta(x) : \theta \in \Theta \} $,
for example, Gaussian processes as in Cakmak et al. \cite{cakmak2024contextual}, stochastic forests as in Kallus and Mao \cite{kallus2023stochastic}, and Gaussian mixture models as in Yoon et al. \cite{yoon2025data}.
For a more detailed discussion of this methodology,
we refer the interested reader to~\cite{bertsimas2020predictive,sadana2025survey}.

The ILO framework has been widely used in operations research and management sciences, yet it remains underexplored from the perspective of optimization theory.
In this paper, we fill out the gap. 
Specifically, 
we derive the first-order necessary optimality conditions 
 in terms of Mordukhovich limiting subdifferentials. To this end, we 
%Under mild assumptions, 
reformulate \eqref{eq:ILO-1} via interchangeability principle \cite{hiai1977integrals} as
\begin{subequations}
\label{eq:ILO-2}
    \begin{eqnarray}
     % (\inmat{ILO--2})\qquad
        \min_{\theta\in\Theta}
        && 
        % \mathbb{E}_{P_X}\left[\mathbb{E}_{P_{Y|X}}[L(z^*(\theta,x),X,Y)]\right] = 
        \mathbb{E}_{P} \left[
        v(\theta,X)\right] \label{eqn:upper-level-problem-B}\\
        \inmat{s.t.} && v(\theta,x) =  \min_{z}\mathbb{E}_{P} [L(z,x,Y,\theta)|X=x]\label{eqn:second-stage-obj-B}\\
        &&\qquad \qquad \inmat{s.t.}\quad z\in\arg \min_{z\in \mathcal{Z}} \mathfrak{C}(z,\theta,x), \; \inmat{for a.e.} \;x\in \mathcal{X}, \quad
        \label{eqn:lower-level-problem-B}
        % x \in \mathcal{X}.
\end{eqnarray}
\end{subequations}
where $v(\theta,x)$ denotes 
the optimal value of the second stage problem \eqref{eqn:second-stage-obj-B}-\eqref{eqn:lower-level-problem-B}. 
In the case where 
$\mathfrak{C}(z,\theta,x)$ is convex and continuously differentiable in $z$, the lower level problem \eqref{eqn:lower-level-problem-B} can be represented by 
%equilibrium 
a variational inequality
constraint derived from the first-order optimality conditions of \eqref{eqn:lower-level-problem-B}.
This yields a two-stage stochastic mathematical program with 
variational inequality
constraints which is analogous 
to 
two-stage stochastic mathematical program with equilibrium constraints (SMPEC), see
\cite{shapiro2008stochastic,xu2010necessary}.
When $\mathfrak{C}(z,\theta,x)$ is nonconvex, 
\eqref{eqn:lower-level-problem-B} cannot 
be equivalently represented by its first-order optimality conditions.
Consequently, we adopt 
a value-function approach proposed by Ye and Zhu \cite{ye1995optimality} to derive an equivalent single-level formulation for the stochastic bilevel problem.
We extend the combined MPEC and value-function approach in Ye and Zhu \cite{ye2010new} to the two-stage stochastic setting.
To obtain nondegenerate multiplier rules for the value-function constraint and hence meaningful first-order optimality conditions, we broaden the classical notion of partial calmness (\cite{mehlitz2021note,ye1995optimality}) to stochastic bilevel optimization by introducing \textit{uniform partial calmness} and \textit{stochastic partial calmness}.
%{ \color{blue}
Building on these techniques, we derive the first-order necessary optimality condition for the ILO problem. 
These conditions can be exploited to develop gradient-based algorithms for solving ILO problems; see, for example, \cite{bai2024alternating,chen2021closing,giovannelli2025inexact,hong2023two,hu2023contextual,liu2020generic}. We leave this for future research.

The main contributions of the paper can be summarized as follows.
\begin{itemize}
         \item[(a)] When the lower-level decision-making problem is convex and continuously differentiable in the decision variable, we reformulate the ILO problem as a two-stage stochastic program with variational inequality constraints \eqref{eq:ILO-3}.
         For the second-stage problem, we use 
         the notion of coderivative multiplier (CD-multiplier) to characterize the stationary system under the well-known no nonzero abnormal multiplier constraint qualification (NNAMCQ).
         Building on sensitivity analysis of the second-stage value function (Proposition \ref{prop:subdifferentiable-V}), we derive a first-order necessary optimality condition for the ILO problem with CD-multipliers (Theorem \ref{thm:first-order-necessary-condition}) that is stated in terms of an upper-level stationary condition involving the expected subgradient of the value function in the sense of Aumann's integral and pointwise M-stationary conditions for second-stage solutions. 
         Moreover, under a uniform inf-compactness condition, we establish the existence of measurable selections of both second-stage solutions and multipliers (Theorem \ref{thm:first-order-necessary-condition-measurable}), which is essential for the asymptotic analysis of sample average approximations.

         \item[(b)] When the lower-level problem is nonconvex, we adopt a value-function approach proposed by Ye and Zhu \cite{ye2010new}  for the ILO problem. 
         To this end, we introduce the stochastic partial calmness condition for the ILO to obtain nondegenerate multiplier rules and discuss its relationship with the classical (uniform) counterpart in a deterministic setting (Proposition \ref{prop:uniform/stochastic-PC}).
         Under this condition, we reformulate the ILO problem as a two-stage stochastic program with variational inequality constraints~\eqref{eqn:CP2-mu-vfunc} and obtain the first-order necessary optimality condition with CD-multipliers by following the scheme in the convex case (Theorem \ref{thm:first-order-necessary-condition2}).
         For the optimality condition, we characterize the subdifferential of the value function and offer a simplified expression of it as a finite convex combination of elements (Proposition \ref{prop:Cstar-lipschitz-differentiable}).
 
     \item [(c)] 
     We apply the established optimality conditions to 
     the “smart predict then optimize” (SPO) model introduced by Elmachtoub and Grigas \cite{elmachtoub2022smart}, which can be formulated as a special case of the ILO model by selecting the specific `SPO' loss in the upper-level objective and has been widely used in operations research and management science.
     In the first illustrative example (Section \ref{section:portfolio}),
     we study the SPO model in a portfolio selection setting, where predicting the conditional distribution of asset returns reduces to predicting their expected values as point estimates.
     We use a linear model to predict these expected returns and derive the first-order necessary optimality condition, in which the Mordukhovich coderivative of the normal cone to the simplex is explicitly characterized in Corollary \ref{corollary:simplex}.
     In the second example (Section \ref{section:newsvendor}), 
     we examine the SPO model in a newsvendor setting and adopt a kernel regression approach to predict the conditional demand distribution.
     We then obtain an optimality condition by invoking the Mordukhovich coderivative of the normal cone to the first orthant (Corollary \ref{corollary:opt-cond-first-orthant}) and present an explicit formula for computing each component when a Gaussian kernel is used.

     \end{itemize}

The rest of this paper is organized as follows. 
In section \ref{sec:preliminary}, we recall some basic notions in variational analysis and nonsmooth analysis 
including coderivatives and Mordukhovich 
subdifferential.
In section \ref{sec:kkt-approach}, we
derive first-order necessary optimality conditions
of ILO problem \eqref{eq:ILO-2} when $\mathfrak{C}(z,\theta,x)$  is convex.
In section \ref{eqn:nonconvex}, 
we extend the discussion to the case that
$\mathfrak{C}(z,\theta,x)$  is nonconvex.
In Section \ref{sec:motivating-applications}, 
 we apply the established optimality conditions to some 
 well-known 
 ILO problems. 

\subsection{Notation}

Throughout this paper, we use the following notation.
For a vector $d\in\mathbb{R}^n$, $\|d\|$ denotes the Euclidean norm of $d$,
$d_i$ is the $i$th component of $d$ and $d_I$ is the subvector composed of the components $d_i,i\in I$. We use $d^\top$ to denote the transpose a vector $d$.
For a $m$-by-$n$ matrix $A$, we denote the $i$th row of $A$ by $a_i^T\in \mathbb{R}^{1\times n}$ with $i=1,\dots, m$ and then $A=(a_1,\dots,a_m)^T$.
For index sets $I\subset \{1,\dots,m\}$, $A_I$ denotes the submatrix of $A$ with rows specified by $I$ and $|I|$ denotes the size of $I$.
$I^{m\times m}$ denotes the identity matrix of size $m$. 
With a slight abuse of notation, $0$ and $1$ may denote all-zero and all-one vectors of appropriate dimensions when the context is clear. Otherwise, we write $\boldsymbol{0}$ and $\boldsymbol{1}$ for such vectors.

For a set-valued mapping $\Phi:\mathbb{R}^m \rightarrow 2^{\mathbb{R}^q}$ (assigning to each $z\in\mathbb{R}^m$ a set $\Phi(z)\subset \mathbb{R}^q$), we denote by $\inmat{gph}\Phi$, i.e., $\inmat{pgh}\Phi:=\{ (z,v)\in \mathbb{R}^m \times \mathbb{R}^q: v\in \Phi(z) \}$.
$\inmat{int}(\mathcal{S})$, $\inmat{cl}(\mathcal{S})$, and $\inmat{co}(\mathcal{S})$ denote the interior, the closure, and the convex hull of a set $\mathcal{S}$, respectively. 
We denote the closed unit ball by $\mathcal{B}$, the small neighborhood of $x\in \mathcal{X}$ relative to $\mathcal{X}$ by $\mathcal{B}(x)$, and the closed ball centered at $(x,z)\in \mathcal{X}\times \mathcal{Z}$ with radius $\delta>0$ by $\mathcal{B}((x,z),\delta)$.
The addition and subtraction of two sets $\mathcal{S}$ and $\mathcal{T}$ are in the sense of Minkowski, 
i.e.,
$\mathcal{S}+\mathcal{T}:= \{s+t: s\in \mathcal{S}, t\in \mathcal{T}\}$ and $\mathcal{S}-\mathcal{T}:= \{s-t: s\in \mathcal{S}, t\in \mathcal{T}\}$.

\section{Mordukhovich Subdifferential}
\label{sec:preliminary}

% {\color{purple}
Since ILO model is fundamentally a bilevel program, 
%the derivation of its
we characterize the first-order optimality conditions
%requires us to adopt 
in terms of Mordukhovich subdifferentials.
We refer interested readers to the monograph \cite{dempe2020bilevel} for a comprehensive overview of bilevel optimization.
In this section, we recall the definitions and related notions in variational analysis and nonsmooth analysis. 
Detailed discussions on these subjects can be found in \cite{clarke1990optimization,clarke1998nonsmooth,mordukhovich2018variational2,rockafellar1998variational}.

\subsection{Variational  analysis}

Let $\Phi: \mathbb{R}^m \rightarrow 2^{\mathbb{R}^m}$ be a set-valued mapping.  Denote by
$
\limsup_{x\rightarrow \bar x} \Phi(x)
$
the {\em  Painlev\'{e}-Kuratowski upper limit}\footnote{In some references, it is also called outer limit, see \cite{rockafellar1998variational}.}
\begin{eqnarray*}
\limsup_{x\rightarrow \bar x} \Phi(x):=  \{ v\in R^m: \exists \mbox{ sequences } x_k\rightarrow  \bar x , v_k \rightarrow v \mbox{ with } v_k  \in \Phi(x_k) \ \  \forall k=1,2,\dots \}.
\end{eqnarray*}
Let $\C$ be a nonempty subset of $\mathbb{R}^m$. Given $z\in$ cl\; $\C$, the convex cone
$$
\N_\C^\pi(z):= \{\zeta \in R^m:
 \exists \sigma >0, \mbox{ such that } \langle \zeta, z'-z\rangle \leq \sigma \|z'-z\|^2 \ \forall z'\in \C\}
 $$
is called the {\em proximal normal cone} to set $\C$ at point $z$, the closed cone
$$
\N_\C(z):= \limsup_{z'\rightarrow  z, z'\in \C}\N_\C^\pi(z')
$$
is called the {\em limiting normal cone} (also known %referred to
as  Mordukhovich normal cone or basic normal cone) to $\C$ at point $ z$.

Note that Mordukhovich originally used  the Fr\'echet  (also called regular) normal cone instead of the proximal normal cone %.
to construct
the limiting normal cone, see \cite[Definition 1.1 (ii)]{mordukhovich2018variational1}.
The two definitions coincide in the finite-dimensional space (see \cite[page 345]{rockafellar1998variational} for a discussion).
The limiting normal cone is in general smaller than the Clarke normal cone,
% {\color{red} HX: cite Clarke's book} 
which is equal to the convex hull $\text{co} \;\N_\C(z)$, see e.g. \cite{clarke1990optimization,clarke1998nonsmooth}. 
In the case when $\C$ is convex, the proximal normal cone, the limiting normal cone and the Clarke normal cone coincide with the normal cone in the sense of the convex analysis, i.e.,
\[
\N_\C(z):=\left \{\zeta \in \bbr^m:\langle \zeta, z'-z\rangle\leq 0,\;\;\forall\,z'\in
\C\right\}.\]
For set-valued mappings, the definition for a limiting normal cone leads to the definition of Mordukhovich coderivative, which is first introduced by Mordukhovich \cite{mordukhovich1980metric}.

\begin{definition}[Coderivatives] \label{def:mordukhovich-coderivative}
Let $\Phi: \mathbb{R}^m \rightarrow 2^{\mathbb{R}^q} $ be an arbitrary set-valued mapping and $(\bar{z},\bar v) \in {\rm cl \; gph } \Phi$. The {\em coderivative} of $\Phi$ at point $(\bar z,\bar v)$ is defined as %The set-valued mapping
$$
D^*\Phi(\bar z,\bar v)(\eta):=\left \{\zeta \in \mathbb{R}^m:  (\zeta, -\eta) \in \N_{gph \Phi}(\bar z,\bar v)\right \}.
$$
 By convention, for $(\bar z,\bar v)\not \in {\rm cl \; gph } \Phi$, %we define
$D^*\Phi(\bar z,\bar v)(\eta)=\emptyset$.

\end{definition}

If $(\bar{z},\bar{v})$ is an interior point of ${\rm gph} \Phi$, then $\N_{gph \Phi}(\bar z,\bar v)=\{(0,0)\}$, which implies $D^*\Phi(\bar z,\bar v)(0)=\{0\}$ and $D^*\Phi(\bar z,\bar v)(\eta)=\emptyset$ for $\eta\neq0$.
A particularly interesting case relevant to our discussions later on is when
$\Phi(z)=\N_\C(z)$  and %where
$\C$ is a closed convex set.
By the definition of coderivative,
$$
\zeta \in D^*\N_\C(\bar z,\bar v) (\eta) \Longleftrightarrow  (\zeta, - \eta) \in \N_{gph \N_\C} (\bar z, \bar v).
$$
Thus
%Hence, the
calculation of the coderivative  $D^*\Phi(\bar z,\bar v)(\eta)$ 
is down to 
the calculation of the limiting normal cone to the normal cone
$ \N_{gph \N_\C} (\bar z, \bar v)$.

\begin{definition}[Subdifferentials]
 Let $f: \mathbb{R}^n \rightarrow \overline \R$ be a lower semicontinuous function and finite at $x\in \mathbb{R}^n$.
 The proximal subdifferential
(\cite[Definition 8.45]{rockafellar1998variational})
 of $f $ at $x$
is defined as
\begin{eqnarray*}
\partial^\pi f(x)&:=&\{\zeta \in \mathbb{R}^n: \exists \sigma>0,\delta >0 \mbox{ such that } f(y) \geq f(x) +\langle \zeta, y-x \rangle -\sigma \|y-x\|^2 \\
&& \qquad \qquad \qquad \forall y \in B(x,\delta)\}
\end{eqnarray*}
and the limiting (Mordukhovich or basic \cite{mordukhovich2018variational1}) subdifferential of $f$ at ${x}$
is defined as
$$
\partial f({x}) :=\limsup_{x'\stackrel{f}{\rightarrow} \ x} \partial^\pi f(x'),
$$
where $x'\stackrel{f}{\rightarrow} \ x$ signifies that $x'$ and $f(x')$ converge to $x$ and $f(x)$ respectively.

When
%the function
$f$ is Lipschitz continuous near ${x}$, the Clarke subdifferential
\cite{clarke1990optimization}
of $f$ at $x$ is equal to
$co \partial f({x})$.
\end{definition}

%30-aug modify the following comments as ref two suggests
Note that Mordukhovich defined the limiting subgradient via Fr\'echet limiting normal cones and
Fr\'echet subgradients (also known as regular subgradients),
see \cite[Theorem 1.89]{mordukhovich2018variational1}. The equivalence of the two definitions is well-known, see the commentary by Rockafellar and Wets \cite[page 345]{rockafellar1998variational}.
The limiting subdifferential is, in general, smaller than the Clarke subdifferential, and in the case when $f$ is convex and
locally  Lipschitz, the proximal subdifferential, the limiting subdifferential, and the Clarke subdifferential
 coincide with the subdifferential in the sense of convex analysis \cite{rockafellar1998variational}. In
the case when $f$ is continuously differentiable, these subdifferentials reduce to
the normal gradient $\nabla f( x)$, i.e., $\partial f(x)=\{ \nabla f( x)\}.$

\subsection{Set-valued mappings and measurability}

Let $\X$ be a closed subset of $\mathbb{R}^n$.
A set-valued mapping $\Phi:\X\to 2^{\mathbb{R}^m}$ is said to be {\em closed}
at $\bar x$ if for $x_k\in\X$, $x_k\to \bar x$, $y_k\in \Phi(x_k)$ and
$y_k\to \bar{y}$ implies $\bar{y}\in \Phi(\bar{x})$.
$\Phi$ is said to be {\em uniformly compact} near $\bar{x}\in\X$ if there is a neighborhood
$B(\bar{x})$ of $\bar{x}$ such that the closure of $\bigcup_{x\in B(\bar{x})} \Phi(x)$
is compact. $\Phi$ is said to be {\em upper semi-continuous} in the sense of Berge \cite{berge1963topological}  at $\bar{x}\in\X$
if for every $\epsilon>0$, there exists a $\delta>0$ such that %neighborhood
%$B(\bar{x})$ of $\bar{x}$ such that
$$
\Phi(\bar{x}+\delta \B) \subset \Phi(\bar{x}) +\epsilon\B,
$$
where $\B$ denotes the closed unit ball in $\mathbb{R}^m$.
% the closure of $\bigcup_{x\in S(\bar{x}} F(x)$
The following result was known, see \cite{gauvin2009differential,hogan1973point}.

\begin{proposition}  Let $\Phi:\X\to 2^{\mathbb{R}^m}$ be uniformly compact near $\bar{x}$. Then
$\Phi$ is upper semi-continuous at $\bar{x}$ if and only if $\Phi$ is closed.
\label{p-usc}
\end{proposition}

Let us now consider %now
a stochastic set-valued mapping.
Let $(\Omega,\F,P)$ be a probability space. For fixed $x$, let $\A(x,\w):\Omega\to 2^{\mathbb{R}^n}$  be a set-valued mapping
%which takes values
whose value is a closed subset of $\mathbb{R}^n$. %Following the notation of Hess \cite{christian2002set}, we
 Let $\mathfrak{B}(\mathbb{R}^n)$ or simply $\mathfrak{B}$  denote the space of closed bounded subsets of $\mathbb{R}^n$
  endowed with topology $\tau_{\mathbb{H}}$ generated by the Hausdorff distance $\mathbb{H}$.
  %{\bf Huifu, do you think that we should define  $\mathbb{H}$?}
  We consider the Borel $\sigma$-field $\G(\mathfrak{B},\tau_{\mathbb{H}})$ generated by
  the $\tau_{\mathbb{H}}$-open subsets of $\mathfrak{B}$. A set-valued mapping  $\A(x,\w):\Omega\to 2^{\mathbb{R}^n}$
  is said to be %{\em strongly
  $\F$-measurable if, for every member $\W$ of $\G(\mathfrak{B},\tau_{\mathbb{H}})$,
  one has $\A^{-1}(\W)\in\F$.
By a {\em measurable selection} of  ${\cal A}(x,\omega)$, we refer to
a vector  $A(x,\omega)\in {\cal A}(x,\omega)$, which
 is measurable. Note that such measurable selections exist if ${\cal A}(x,\omega)$ is measurable,
see \cite{artstein1975strong} and references therein. %\cite{aum65,kr65}.

For a general set-valued mapping which is not necessarily measurable,
the {\em expectation of ${\cal A}(x,\omega)$}, denoted by
$\bbe[{\cal A}(x,\omega)]$, is defined as
the collection of $\bbe[A(x,\omega)]$ where $A(x,\omega)$
is an integrable selection, and the integrability is in the sense of Aumann \cite{aumann1965integrals}.
$\bbe[{\cal A}(x,\omega)]$ is regarded as
{\em well-defined}
if
% $\bbe[{\cal A}(x,\omega)]\in \mathfrak{B}$
it is
%a
nonempty.
%30-aug add
%and compact
% set of $\mathbb{R}^n$.
A sufficient condition of the well-definedness of the expectation
is that ${\cal A}(x,\omega)$ is measurable and
$
\bbe[\|{\cal A}(x,\omega)\|] :=
\bbe[\mathbb{H}(0, {\cal A}(x,\omega)]<\infty,
$
in which case  $\bbe[{\cal A}(x,\omega)]\in \mathfrak{B}$.
See \cite[Theorem 2]{aumann1965integrals}. In such a case, $\A$ is called {\em integrably
bounded} in \cite{aumann1965integrals,christian2002set}.
The following result is well-known, see e.g.
%30-aug address ref 2's comment
 \cite[Page 307, Section 8.1]{aubin1992set} and
 \cite[Lemmas 3.1-3.2]{christian2002set}.

\section{First-order Necessary Optimality Condition: 
Convex Case}
\label{sec:kkt-approach}

In this section, we derive the first-order necessary optimality condition of problem \eqref{eq:ILO-2}. We concentrate on the case  where the lower-level decision-making problem \eqref{eqn:lower-level-problem-B} is convex and  
discuss the nonconvex case in Section 4 as it requires 
different mathematical treatments.  
To this end, we need to make the following assumption.

\begin{assumption} \label{assumption:convex-smooth}
    For any $(\theta,x)\in \Theta \times \mathcal{X}$, $\mathfrak{C}(z,\theta,x)$ is convex and continuously differentiable in $z$. The feasible set $\mathcal{Z}$ is a closed, convex subset of $\mathbb{R}^{d_z}$.
\end{assumption}
Under Assumption~\ref{assumption:convex-smooth}, 
the ILO problem \eqref{eq:ILO-2} can be equivalently formulated via the first-order optimality condition of the lower level decision-making problem as 
\begin{subequations}
\label{eq:ILO-3}
    \begin{eqnarray}
     % (\inmat{ILO--convex})\qquad
        \min_{\theta\in\Theta}
        && 
        % \mathbb{E}_{P_X}\left[\mathbb{E}_{P_{Y|X}}[L(z^*(\theta,x),X,Y)]\right] = 
        \mathbb{E}_{P} \left[
        v(\theta,X)\right] \label{eqn:upper-level-problem-C}\\
        \inmat{s.t.} && v(\theta,x) :=  \min_{z}\mathbb{E}_{P} [L(z,x,Y,\theta)|X=x]\label{eqn:second-stage-obj-C}\\
        &&\qquad \qquad \inmat{s.t.}\quad 0\in \nabla_z\mathfrak{C}(z,\theta,x) + \mathcal{N}_\mathcal{Z} (z),
        % z\in\arg        \min_{z\in \mathcal{Z}} \mathfrak{C}(z,\theta,x), 
        \; \inmat{for a.e.} \;x\in \mathcal{X}. \quad
        \label{eqn:lower-level-optimal-C}
        % x \in \mathcal{X}.
\end{eqnarray}
\end{subequations}
 Let $\Gamma(\theta,x)$ be the set of optimal solutions of problem \eqref{eqn:second-stage-obj-C}–\eqref{eqn:lower-level-optimal-C}.
 %at $(\theta,x)$.
We begin by introducing the coderivative multiplier of 
the lower-level problem \eqref{eqn:second-stage-obj-C}-\eqref{eqn:lower-level-optimal-C}.

\begin{definition}[M-multipliers] \label{def:cd-multiplier}
Let $(\theta,x)\in \Theta\times \mathcal{X}$ be fixed and  $z$ be a
feasible solution of the second stage problem \eqref{eqn:second-stage-obj-C}-\eqref{eqn:lower-level-optimal-C}.
 $z$ is said to be a Mordukhovich stationary point, and $\eta\in \mathbb{R}^{d_z}$ is an M-multiplier of the second stage problem \eqref{eqn:second-stage-obj-C}-\eqref{eqn:lower-level-optimal-C} at $z$ if they satisfy the following equation
\begin{eqnarray}
0& \in& \nabla_{z} \mathbb{E}_P\left[L(z,X,Y,\theta)\mid X=x\right] + \nabla_{zz}^2 \mathfrak{C}(z(x),\theta,x)^{\top} \eta  + D^*\N_\mathcal{Z}\big({z}, - \nabla_{z}\mathfrak{C}(z,\theta,x)\big) (\eta), \qquad
\end{eqnarray}
where $D^*$ is defined in Definition \ref{def:mordukhovich-coderivative}.
Let ${\cal M}(\theta, x, z)$  denote the set of M-multipliers at the stationary point $z$ for the given $(\theta,x)$.
\end{definition}
Note that the coderivative multiplier 
%(CD-multiplier in \cite{lucet2001sensitivity}, 
is also known as  M-multiplier, see e.g.~\cite{xu2010necessary}.
%is defined as follows.
To ensure non-emptiness of 
set ${\cal M}(\theta, x, z)$ at the local optimal solution $z$ of the second-stage problem,
we adopt the well-known
no nonzero abnormal multiplier constraint qualification (NNAMCQ).

\begin{definition}[NNAMCQ]
    %We say that
    Problem \eqref{eqn:second-stage-obj-C}-\eqref{eqn:lower-level-optimal-C} is said to satisfy 
    NNAMCQ at a feasible point $z$ of  if
    \begin{eqnarray*}
        0\in \nabla_{zz}^2 \mathfrak{C}(z(x),\theta,x)^{\top} \eta  + D^*\N_\mathcal{Z}({z}, - \nabla_{z}\mathfrak{C}(z,\theta,x)) (\eta) \Longrightarrow \eta=0. 
    \end{eqnarray*}
\end{definition}
The notion of NNAMCQ is introduced by Ye in \cite{ye1997necessary} and is used 
to derive
%in derivation of 
the first-order necessary optimality condition of two-stage
SMPECs by Xu and Ye in 
\cite{xu2010necessary}.
%,ye2000constraint,ye1997necessary},
The constraint qualification coincides with coderivative nonsingularity in \cite{rockafellar1998variational} 
and is equivalent to the metric regularity constraint qualification; see, e.g., \cite[Theorem 9.43]{rockafellar1998variational} or \cite[Theorem 4.2]{dontchev2014implicit}. 
%In practice, 
The NNAMCQ is guaranteed by 
%(and metric regularity) can often be verified via many classical and 
some verifiable sufficient conditions. 
For instance, when $\nabla_{z} \mathfrak{C}(z(x),\theta,x) + \mathcal{N}_\mathcal{Z}(z)$ is single-valued and continuously differentiable, NNAMCQ reduces to the full rank of the Jacobian at $z$. 
In the case that the feasible set is characterized by a system of inequalities,  NNAMCQ is sufficed by 
classical constraint qualifications in nonlinear programming such as MFCQ/LICQ.
%, provide convenient sufficient conditions for metric regularity. 
Moreover, the second order sufficient condition (SOSC) and LICQ imply
strong regularity in the sense of Robinson \cite[Theorem 4.1]{robinson1980strongly} and hence NNAMCQ, see, e.g., \cite[Theorem 4.7]{ye2000constraint}, \cite[Theorem 3.2]{ye1997necessary}, \cite[Proposition 2.13]{xu2010necessary}, or monographs \cite[Section 13.3]{bonnans2006numerical}, \cite[Section 9F,9G]{rockafellar1998variational}, \cite[Chapter 4]{dontchev2014implicit}.

\begin{assumption} \label{assumption:differentiable-regular}
   Let  
$
\mathcal{Z}^*(\theta,x) := 
\arg \min_{z\in \mathcal{Z}} \mathfrak{C}(z,\theta,x)
$
for $(\theta,x)\in \Theta \times \mathcal{X}$.
For $\theta\in \Theta$ and a.e. $x\in \mathcal{X}$, the NNAMCQ holds at $z^*(\theta,x)\in \mathcal{Z}^*(\theta,x)$.
\end{assumption}

Under Assumption \ref{assumption:convex-smooth}, $\mathcal{Z}^*(\theta,x)= \left\{z\in \mathcal{Z}: 0\in \nabla_z \mathfrak{C}(z,\theta,x) +{\cal N}_{\cal Z}(z)\right\}$ for $(\theta,x) \in \Theta \times \mathcal{X}$.
Next, we make assumptions on the continuous differentiability 
of 
 $\nabla_z \mathfrak{C}(z,\theta,x)$ and $\mathbb{E}_P[L(z,X,Y,\theta)|X=x]$ to ensure
 the M-subdifferential of the value function $v(\theta,x)$ with respect to $\theta$ is bounded.

\begin{assumption} \label{assumption:lipschitz_in_z}
    For any $x\in\mathcal{X}$,
    % (a) 
    $\nabla_z \mathfrak{C}(z,\theta,x)$ and $\mathbb{E}_P[L(z,X,Y,\theta)|X=x]$ are continuously differentiable over $\mathcal{Z}\times \Theta$.
    % and locally Lipschitz continuous in $(z,\theta)$.
\end{assumption}

%In fact, t
It is possible to relax the condition to
local Lipschitz continuity 
in which case the optimality condition can be derived 
by virtue of \cite[Theorem 3.6 and Corollary 3.7]{lucet2001sensitivity}.
To provide a clearer characterization of this requirement in applications, the next proposition gives a sufficient condition ensuring local Lipschitz continuity with respect to $\theta$ in the setting where $\nabla_z \mathfrak{C}(z,\theta,x):=\mathbb{E}_{f_\theta(x)}[\nabla_z c(z, Y)]$.

\begin{proposition} \label{prop:lipschitz}
Assume:
(a) For any fixed $(x,y)\in \mathcal{X}\times \mathcal{Y}$, the probability density function of $f_\theta(x)$, denoted by $p_\theta(y;x)$, is continuously differentiable in $\theta$; 
    (b) for any $(\tilde{\theta},\tilde{z})\in \Theta\times \mathcal{Z}$, there exists an integrable function $g(y;x)$, i.e. $L_g= \int_\mathcal{Y} g(y;x)dy <\infty$, such that for  $\theta$ in a neighborhood of $\tilde{\theta}$ and $z$ in a neighborhood of $\tilde{z}$,
    \begin{eqnarray}
        \left\| \nabla_z c(z,y) \right\| \left\| \nabla_\theta p_{\theta}(y;x) \right\| \leq g(y;x), \quad \forall (x,y)\in \mathcal{X}\times \mathcal{Y};
    \end{eqnarray}
   (c) Assumption \ref{assumption:convex-smooth} holds.
   Then
    %and Assumption \ref{assumption:continuous-density}, 
    $\nabla_z \mathfrak{C}(z,\theta,x)=\mathbb{E}_{f_\theta(x)}[ \nabla_z c(z,Y)]$ is locally Lipschitz continuous in $\theta\in \Theta$.
\end{proposition}
Note that if $\left\| \nabla_z c(z,y) \right\|$ is bounded uniformly
%bounded 
for all $y\in \mathcal{Y}$ and $z\in \mathcal{Z}$, condition (b) of Proposition \ref{prop:lipschitz} reduces to the requirement of $\int_\mathcal{Y} \left\| \nabla_\theta p_{\theta}(y;x) \right\| dy < \infty$,
where $\int_\mathcal{Y} \| \nabla_\theta p_\theta(y) \| dy$ can be rewritten as  $\mathbb{E}_{f_\theta(x)}[\| \nabla_\theta \log p_\theta(y) \|]$. By %applying 
Jensen's inequality,
%we have 
\begin{eqnarray*}
    \mathbb{E}_{f_\theta(x)}\left[ \| \nabla_\theta \log p_\theta(y) \| \right] \le \sqrt{ \mathbb{E}_{f_\theta(x)} \left[ \| \nabla_\theta \log p_\theta(y) \|^2 \right] } = \sqrt{\text{Trace}(I({f_\theta(x)}))}
\end{eqnarray*}
where $I(f_\theta(x))$ denotes the Fisher information matrix of $f_\theta(x)$ 
(\cite[equation (2.6.8)]{lehmann1998theory})
defined by $$I_{ij}({f_\theta(x)}) = \mathbb{E}_{f_\theta(x)}\left[ \nabla_{\theta_i} \log p_\theta(Y;x) \nabla_{\theta_j} \log p_\theta(Y;x) \right].$$ 
Consequently, the uniform boundedness assumption holds whenever the Fisher information is finite, which is commonly satisfied, e.g., for distributions within the exponential family \cite[Section 2.5]{lehmann1998theory}.

\noindent
\textbf{Proof.}
For $\theta_1,\theta_2\in \Theta$ in a neighborhood of $\tilde{\theta}$, 
\begin{eqnarray*}
    \left\|\nabla_z\mathbb{E}_{f_{\theta_1}(x)}\left[c(z,Y)\right] - \nabla_z\mathbb{E}_{f_{\theta_2}(x)}\left[c(z,Y)\right]\right\|
    & = &
    \left\|\mathbb{E}_{f_{\theta_1}(x)}\left[\nabla_zc(z,Y)\right] - \mathbb{E}_{f_{\theta_2}(x)}\left[\nabla_zc(z,Y)\right]\right\| \\
    &=& \left\| \int_\mathcal{Y} \nabla_z c(z,y)\left( p_{\theta_1}(y;x) - p_{\theta_2}(y;x) \right) dy \right\|\\
    &\leq & \int_\mathcal{Y} \left\|\nabla_z c(z,y) \right\| \left\| p_{\theta_1}(y;x) - p_{\theta_2}(y;x) \right\| dy\\
    &\leq & \int_\mathcal{Y} \left\|\nabla_z c(z,y) \right\| \left\| \nabla_\theta p_{\tilde{\theta}}(y;x) \right\| \|\theta_1-\theta_2\| dy\\
    &\leq & \int_\mathcal{Y}g(y;x) dy  \|\theta_1-\theta_2\| = L_g\|\theta_1-\theta_2\|.
\end{eqnarray*}
The first equation is derived from Assumption \ref{assumption:convex-smooth}, in which the interchangeability between expectation and derivative can be ensured by the dominated convergence theorem.
The second inequality is derived using the mean value theorem with $\tilde{\theta}\in [\theta_1,\theta_2]$, and the last inequality follows from condition (b).
This result indicates that $\nabla_z\mathbb{E}_{f_{\theta_1}(x)}\left[c(z,Y)\right]$ is locally Lipschitz continuous in $\theta$. 
\hfill $\Box$

To ensure the existence of an optimal solution 
and the boundedness of 
the set of optimal solutions
%set for 
of the second-stage problem \eqref{eqn:second-stage-obj-C}-\eqref{eqn:lower-level-optimal-C}, we 
%assume 
need the following inf-compact conditions additionally.

\begin{assumption} \label{assumption:inf-compactness}
    Let $\theta \in \Theta$ and $x\in \mathcal{X}$ be fixed. There exists a constant $\delta>0$ such that the set 
    \begin{eqnarray*}
        \Big\{ z: r\in  \nabla_z \mathfrak{C}(z,\theta,x) +{\cal N}_{\cal Z}(z), \mathbb{E}_{P} [L(z,x,Y,\theta)|X=x] \leq \alpha, r\in \delta\mathcal{B} \Big\}
    \end{eqnarray*}
    is bounded for some constant $\alpha$.
\end{assumption}

Based on the 
%above
assumption and previous assumptions, we are ready to 
address the local Lipschitz continuity and subdifferentiability 
of the value function $v(\theta,x)$ in the next proposition.

\begin{proposition} \label{prop:subdifferentiable-V}
    Let $\bar{\theta}\in \Theta$ be fixed. 
    Suppose (a) Assumptions \ref{assumption:convex-smooth}-\ref{assumption:lipschitz_in_z} hold; (b) Assumption \ref{assumption:inf-compactness} holds at $(\bar{\theta},x)$ for every $x\in \mathcal{X}$. 
    Then the following assertions hold.
    \begin{itemize}
        \item[(i)] 
        The optimal value mapping $v(\theta,x)$ is locally Lipschitz continuous over $\Theta\times \mathcal{X}$,
        % near $(\bar\theta,\bar x)$ for any $\bar x\in \mathcal{X}$, 
        and $v(\bar{\theta},\cdot):\mathcal{X}\rightarrow  \mathbb{R}$ is measurable.

        \item[(ii)] 
        The optimal solution mapping $\Gamma(\theta,\cdot): \mathcal{X}\rightarrow 2^{\mathbb{R}^{d_z}}$ is measurable and compact-valued.
        
        % {\color{red} Obvious? 
        % Is it a closed set-valued mapping?}
    
        \item[(iii)] $v(\theta,x)$ is lower semicontinuous near $\bar \theta$, $\partial_\theta v(\bar \theta,x) \subset \Psi(\bar \theta,x)$ for every $x$, 
        where
        \begin{eqnarray*}
        \Psi(\theta,x):= 
        \bigcup_{z\in\Gamma(\theta,x)} 
        \bigcup_{\eta\in {\cal M}(\theta,x,z)} \left\{
        \nabla_\theta \mathbb{E}_P\left[L(z,X,Y,\theta)\mid X=x\right] + 
        \nabla_{z\theta}^2 \mathfrak{C}(z,\theta,x)^{\top} \eta\right\},
        \end{eqnarray*}
        % \item[(iv)]
        %Then 
        % $v(\theta,\cdot):\mathcal{X}\rightarrow \mathbb{R}$, 
        and $\partial_\theta v(\theta,\cdot):\mathcal{X}\rightarrow 2^{\mathbb{R}^{d_\theta}}$ 
        is measurable.
        
        \item[(iv)]  Let $V(\theta):=\mathbb{E}_P[v(\theta,X)]$. 
        % {\color{blue}
% If $v(\theta,x)$ is Lipschitz in $\theta$ for each fixed $x$ and}
        %$V(\theta)$ is well-defined and 
        Suppose that
        the Lipschitz modulus of $v(\theta,x)$ at $\bar \theta$
        %of $v(\theta,x)$ in $\theta$ 
        is bounded by an integrable function $\kappa(x)$, then $V(\theta)$ is locally Lipschitz continuous at $\bar\theta$ and 
        \begin{eqnarray}
            \partial V(\theta) \subseteq \mathbb{E}_P[\partial_\theta v(\theta,x)].
        \end{eqnarray}
        Moreover, if $v(\theta,x)$ is 
        Clarke regular (Definition 2.3.4 in \cite{clarke1990optimization})
% {\color{red} Clarke regular Definition??} 
in $\theta$, then $V(\theta)$ is Clarke regular and the equality holds.
    \end{itemize}
\end{proposition}

%{\color{red} HX: please provide a proof}

\noindent
\textbf{Proof.}
The proof is similar to that of \cite[Corollary 3.8]{lucet2001sensitivity} and \cite[Proposition 2.15, Theorem 2.16]{xu2010necessary},  we skip the details.
%we arrive at the following results.
\hfill $\Box$

Combined with the definition of CD-multiplier, we summarize the necessary first-order optimality condition for the ILO model.

\begin{theorem} \label{thm:first-order-necessary-condition}
    Assume the conditions of Proposition \ref{prop:subdifferentiable-V} are satisfied.
    Assume also that for fixed $\theta$, there exists a nonnegative function $\sigma(x)$ with $\mathbb{E}_P[\sigma(X)]< \infty$ such that 
    \begin{eqnarray} \label{eqn:theorem-lipshcitz-bound1}
        % \max\left\{ \| \nabla_\theta\mathbb{E}_P\left[L(z,X,Y)\mid X=x\right] \| , 
        \| \nabla_{z\theta}^2 \mathfrak{C}(z,\theta,x) \| 
        % \right\} 
        \leq \sigma(x)
    \end{eqnarray}
    uniformly for all $z\in \Gamma(\theta,x)$.
    If $(\theta,z(\cdot))$ is a local optimal solution of problem \eqref{eq:ILO-3}, then
    % Then for a local optimal solution $\theta^*$ of the ILO model \eqref{eq:ILO-3}, 
    there exist 
    % a selection of stationary point $z(x)$ of problem \eqref{eqn:second-stage-obj-C}-\eqref{eqn:lower-level-optimal-C} and the corresponding 
    a selection of CD multipliers $\eta(x)$ such that 
    \begin{eqnarray}
    \label{eqn:first-order-necessary-condition-convex}
        \begin{cases}
            0 \in  \mathbb{E}_P \left[ 
            \nabla_\theta \mathbb{E}_P\left[L(z(X),X,Y,\theta)\mid X\right] + 
            \nabla_{z\theta}^2 \mathfrak{C}(z(X),\theta,X)^{\top} \eta(X) \right] + \mathcal{N}_\Theta(\theta),\\
            0 \in \nabla_{z} \mathbb{E}_P\left[L(z(X),X,Y,\theta)\mid X=x\right] + \nabla_{zz}^2 \mathfrak{C}(z(x),\theta,x)^{\top} \eta(x)  \\
            \qquad+ D^*\N_\mathcal{Z}(z(x), - \nabla_{z}\mathfrak{C}(z(x),\theta,x)) (\eta(x)),  \; \forall  x\in \mathcal{X},\\
            0 \in \nabla_z \mathfrak{C}(z(x),\theta,x) +{\cal N}_{\cal Z}(z(x)), \quad  \forall x\in \mathcal{X}.
        \end{cases}
    \end{eqnarray}
\end{theorem}

\noindent
\textbf{Proof.}
By Proposition \ref{prop:subdifferentiable-V}(ii) and Assumption \ref{assumption:differentiable-regular}, the set of multiplier $\bigcup_{z\in \Gamma(\theta,x)}M(\theta,x,z)$ is nonempty and compact (\cite[Proposition 2.12]{xu2010necessary}).
This result and condition \eqref{eqn:theorem-lipshcitz-bound1} ensure the Lipschitz continuity of $v(\theta,x)$ in $\theta$ and thus guarantee the assumption in Proposition \ref{prop:subdifferentiable-V}(iv).
Then the result immediately follows from Proposition \ref{prop:subdifferentiable-V}(iv) and the definition of CD-multipliers in Definition \ref{def:cd-multiplier}.
\hfill $\Box$

It is worth noting that in Theorem \ref{thm:first-order-necessary-condition}, we do not require the measurability of the selection of the CD multiplier $\eta(\cdot)$ and thus the well-definedness of 
\begin{eqnarray*}
    \mathbb{E}_P \left[ 
            \nabla_\theta \mathbb{E}_P\left[L(z(X),X,Y,\theta)\mid X\right] + 
            \nabla_{z\theta}^2 \mathfrak{C}(z(X),\theta,X)^{\top} \eta(X) \right].
\end{eqnarray*}
% $\mathbb{E}_P \left[\nabla_{z\theta}^2 \mathfrak{C}(z(X),\theta,X)^{\top} \eta(X) \right]$.
When the underlying distribution $P$ is finitely supported, such as in the case of an empirical distribution within the sample average approximation (SAA)
% {\color{red}Stochastic Approximation Algorithm (SAA) 
% HX: SAA usually refers to sample average approximation
% }
method, the expectation operator simplifies to a finite weighted sum, which allows the selection of multipliers as a multiplier vector in each scenario.
In this case, the measurability of the selection of the CD multipliers holds automatically, and Theorem \ref{thm:first-order-necessary-condition} suffices for formulating the necessary first-order optimality conditions.
However, for continuous underlying distributions and for the asymptotic analysis of the SAA method, it becomes essential to ensure that the second-stage solutions $z(x)$ and multipliers $\eta(x)$ can be chosen as measurable selections. 
In the following discussion, we impose the uniform inf-compactness assumption and ensure the measurability of the selections.

\begin{assumption} [Uniform inf-compactness] \label{assumption:uni-inf-compact}
    Let $\theta \in \Theta$ be fixed. For every $x\in \mathcal{X}$, there exists a constant $\delta>0$ such that the set 
    \begin{eqnarray}
        \Big\{ z: r\in  \nabla_z \mathfrak{C}(z,\theta,x') +{\cal N}_{\cal Z}(z), \mathbb{E}_{P} [L(z,x',Y,\theta)|X=x'] \leq \alpha, r\in \delta\mathcal{B}
        % \mathbb{B}(0,\delta) 
        \Big\}
    \end{eqnarray}
    is bounded for some constant $\alpha$ and every $x'$ in a close neighborhood of $x$ relative to $\mathcal{X}$.
\end{assumption}

Based on this assumption, we establish the upper semicontinuity of $M(\theta,\cdot,\cdot)$ and then obtain a stronger version of Theorem \ref{thm:first-order-necessary-condition} with a measurable multiplier selection $\eta(\cdot)$. 
For $x\in \mathcal{X}$, let 
\begin{eqnarray*}
    \mathcal{H}(x):= \{ \{x'\}\times\Gamma(\theta,x'): x'\in \mathcal{B}(x)\}
\end{eqnarray*}
be a collection of sets in space $\mathbb{R}^{d_z}\times \mathcal{X}$.
Then for $z\in \Gamma(\theta,x)$, we say $M(\theta,\cdot,\cdot)$ is upper semicontinuous at $(z,x)$ relative to set $\mathcal{H}$ if for every $\nu>0$, there exists $\delta>0$ such that
\begin{eqnarray*}
    M(\theta,x',z')\subset M(\theta,x,z)+ \nu\mathcal{B}
\end{eqnarray*}
for all $(x',z')\in \inmat{cl} \mathcal{B}((x,z),\delta)\cap \mathcal{H}$, where $\mathcal{B}$ denotes the closed unit ball and $\mathcal{B}((x,z),\delta)$ denotes a closed ball with radius $\delta$ and center $(x,z)$.

\begin{theorem} \label{thm:first-order-necessary-condition-measurable}
    Let $(\theta,z(\cdot))$ be a local optimal solution of problem \eqref{eq:ILO-3}. Suppose that the conditions of Theorem \ref{thm:first-order-necessary-condition} are satisfied and Assumption \ref{assumption:uni-inf-compact} holds at $\theta$.
    Then the following assertions hold.
    \begin{itemize}
    \item[(i)] Denote the set of multipliers at $(\theta,x)$ by $\mathcal{M}(\theta,x):= \bigcup_{z\in \Gamma(\theta,x)} M(\theta,x,z)$. Then $\mathcal{M}(\theta,x)$ is locally bounded on $\mathcal{X}$, i.e., for any $x\in \mathcal{X}$, $\bigcup_{x'\in \mathcal{B}(x)}\mathcal{M}(\theta,x')$ is bounded.
    
        \item[(ii)] 
        For 
        %Let
        every $x\in\mathcal X$ and $z\in \Gamma(\theta,x)$,
        %. Then 
        $M(\theta,\cdot,\cdot)$ is upper semicontinuous at $(x,z)$ relative to $\mathcal{H}$.

        \item[(iii)] $\Psi(\theta, X)$ is measurable and integrably bounded.

        \item[(iv)] There exists measurable selections $z(\cdot)$ with $z(x)\in \Gamma(\theta,x)$ and $\eta(\cdot)$ with $\eta(x)\in M(\theta,x,z(x))$ such that \eqref{eqn:first-order-necessary-condition-convex} holds.
    \end{itemize}
\end{theorem}

\noindent
\textbf{Proof.}
{\color{black}
Part (i).
Suppose, for the sake of 
a contradiction, that the set of 
multiplier 
%sets 
$\mathcal{M}(\theta,x)$ 
%are
is not locally bounded 
{\color{black} near $\bar x\in \mathcal{X}$}.
% {\color{red} the solution graph??}.
Then there exist sequences $\{x_k\}$, $\{z_k\}$, and $\{\eta_k\}$,  such that $x_k\rightarrow \bar x$, $z_k\in\Gamma(\theta,x_k)$ and $\eta_k\in M(\theta,x_k,z_k)$
with $\|\eta_k\|\to\infty$ and $(x_k,z_k,\eta_k)$ satisfies
%By Definition \ref{def:cd-multiplier}, each pair $(z_k,\eta_k)$ satisfies the following inclusion
\begin{eqnarray}  \label{eqn:proof-measurability}
    0\in \nabla_z \mathbb{E}_P\left[L(z_k,X,Y,\theta)\mid X=x_k\right]
+ \nabla^2_{zz}\mathfrak{C}(z_k,\theta,x_k)^{\top}\eta_k
+ D^{*}\mathcal{N}_\mathcal{Z}\!\left(z_k,-\nabla_z \mathfrak{C}(z_k,\theta,x_k)\right)(\eta_k).
\end{eqnarray}
Let $\tilde{\eta}_k := \eta_k/\|\eta_k\|$.
By the compactness of $\bigcup_{x'\in \mathcal{B}(x)}\Gamma(\theta,x')$ under Assumption \ref{assumption:uni-inf-compact} and the boundedness of $\{\tilde{\eta}_k\}$ due to $\left\|\tilde\eta_k\right\|=1$, 
% there exist subsequences $\{z_{k_j}\}$ and $\{\tilde{\eta}_{k_j}\}$ such that 
we assume for simplicity of exposition that $z_{k} \rightarrow \bar{z}$ and $\tilde{\eta}_{k}\rightarrow \bar{\eta}$ with $\|\bar{\eta}\|=1$ by taking subsequences if necessary.
Since $z_k\in\Gamma(\theta,x_k)$ and %closedness of 
$\Gamma(\theta,\cdot)$ is closed by Proposition \ref{p-usc} and Proposition \ref{prop:subdifferentiable-V} (ii), then 
%we have
$\bar z\in\Gamma(\theta,x)$.
Dividing both sides of \eqref{eqn:proof-measurability}
by $\|\eta_k\|$ and driving $k\rightarrow \infty$, we
%Then we 
obtain a nonzero $\bar\eta\neq 0$ such that
\begin{eqnarray*}
    0\in \nabla^2_{zz}\mathfrak{C}(\bar z,\theta,x)^{\top}\bar\eta
+ D^{*}\mathcal{N}_\mathcal{Z}\left(\bar z,-\nabla_z \mathfrak{C}(\bar z,\theta,x)\right)(\bar\eta),
\end{eqnarray*}
which contradicts NNAMCQ in Assumption \ref{assumption:differentiable-regular} at $\bar{z}$.
This implies the local boundedness of $\mathcal{M}(\theta,x)$.
}

Part (ii).
Let $\{x_k\} \subset \mathcal{B}(x)$ be a sequence 
such that $x_k\rightarrow x$ as $k\rightarrow \infty$.
Let $\{z_k\}$ and $\{\eta_k\}$ such that $z_k\in\Gamma(\theta,x_k)$ and $\eta_k\in M(\theta,x_k,z_k)$.
We show that, 
% {\color{red} $\{\eta_k\}$ is relatively compact??} 
% and 
 for any cluster point $\bar \eta$ of $\{\eta_k\}$, there exists $\bar z\in \mathcal{Z}$ such that $\bar \eta$ belongs to $M(\theta,x,\bar z)$.
By Assumption \ref{assumption:uni-inf-compact} and Part (i), 
$\{z_k\}$ and $\{\eta_k\}$ are all bounded.
Taking a subsequence if necessary, 
we assume for the simplicity of exposition that 
$z_k\rightarrow \bar z$ and $\eta_k\rightarrow \bar\eta$.
By Definition \ref{def:cd-multiplier} that $(z_k,\eta_k)$ satisfies the inclusion \eqref{eqn:proof-measurability},
we drive $k\to\infty$ in \eqref{eqn:proof-measurability} and then obtain
\begin{eqnarray*} 
    0\in \nabla_z \mathbb{E}_P\left[L(\bar z,X,Y,\theta)\mid X=x\right]
+ \nabla^2_{zz}\mathfrak{C}(\bar z,\theta,x)^{\top}\bar \eta
+ D^{*}\mathcal{N}_\mathcal{Z}\!\left(z_k,-\nabla_z \mathfrak{C}(\bar z,\theta,x)\right)(\bar \eta),
\end{eqnarray*}
which means
$\bar\eta\in M(\theta,x,\bar z)$.
In addition,
$\bar z\in\Gamma(\theta,x)$ as proved in part (i). 
Therefore, the graph of $M(\theta,\cdot,\cdot)$ restricted to $\mathcal{H}(x)$ is closed and
uniformly compact near $(x,z)$. 
By Proposition \ref{p-usc}, these properties give rise to the upper semicontinuity
of $M(\theta,\cdot,\cdot)$ at $(x,z)$ relative to $\mathcal{H}$.

Part (iii). 
% To ensure %the 
% measurability of $\Psi(\theta,\cdot)$ {\color{blue} in $x$}, 
We first demonstrate the upper semicontinuity of the set-valued mapping $\Psi(\theta,\cdot)$.
To invoke Proposition \ref{p-usc} for this result, we need to establish the locally uniform compactness and closedness of $\Psi(\theta,\cdot)$.
In part (i), we have already established the 
% This implies the 
local boundedness of $\mathcal{M}(\theta,x)$.
Combining this result with the compactness of $\bigcup_{x'\in \mathcal{B}(x)}\Gamma(\theta,x')$, we establish the boundedness of $\bigcup_{x'\in B(x)}\Psi(\theta,x')$.

To demonstrate the closedness of $\Psi(\theta,\cdot)$, 
we let $\{x_k\} \subset \mathcal{B}(x)$ 
be such that $x_k\rightarrow x$ as $k\rightarrow \infty$ and take sequences $\{p_k\}$ such that $p_k\in\Psi(\theta,x_k)$.
By definition, for each $k$ there exist $z_k\in\Gamma(\theta,x_k)$ and $\eta_k\in M(\theta,x_k,z_k)$ such that $p_k=\nabla_\theta \mathbb{E}[L(z_k,X,Y,\theta)|X=x_k] +\nabla^2_{z\theta}C(z_k,\theta,x_k)^{\top}\eta_k$. Then it suffices to show the closedness of the set of multipliers $\bigcup_{x'\in \mathcal{B}(x)}\mathcal{M}(\theta,x')$.
Since $\{z_k\}$ and $\{\eta_k\}$ are bounded as aforementioned, $z_k\to \bar z$ and $\eta_k\to\bar\eta$ by taking a subsequence if necessary.
Substituting them into \eqref{eqn:proof-measurability} and taking the limit on both sides of the equation, we can see $\bar{\eta}\in M(\theta,\bar{z},x)\subset \mathcal{M}(\theta,x)$ and hence the closedness.

By Proposition \ref{p-usc}, these results lead to the upper semicontinuity of the set-valued mapping $\Psi(\theta,\cdot)$ at $x$.
%Then t
The measurability follows directly from \cite[Corollary 14.14]{rockafellar1998variational} by considering $\Psi(\theta,X)$ as a composition of an upper semicontinuous set-valued mapping $\Psi(\theta,\cdot)$ and a measurable random vector $X$.

Part (iv). 
By the definition of Aumann's integral and Proposition~\ref{prop:subdifferentiable-V} (iv), there exists a measurable selection $q(\cdot)$ such that $q(x)\in \partial_\theta v (\theta,x) \subset \Psi(\theta,x)$ satisfies $\mathbb{E}_P[q(X)]\in \partial_\theta V(\theta)$.
Therefore, for a.e. $x$,
\begin{equation*}
q(x)\in \bigcup_{z\in\Gamma(\theta,x)} \Big(\nabla_\theta \mathbb{E}_P\left[L(z,x,Y,\theta)\mid X=x\right] +\nabla^2_{z\theta}C(z,\theta,x)^{\top} M(\theta,x,z)\Big).
\end{equation*}
We can consider the above inclusion as 
\begin{eqnarray*}
    q(x)\in R\left(x,\mathcal U(x)\right),
\end{eqnarray*}
where $R(x,(z,\eta)):= \nabla_\theta \mathbb{E}_P\left[L(z,x,Y,\theta)\mid X=x\right] +\nabla^2_{z\theta}C(z,\theta,x)^{\top}\eta$ and
$\mathcal U(x):=\{(z,\eta)\in Z\times\mathbb{R}^{d_z}: z\in\Gamma(\theta,x), \eta\in M(\theta,x,z)\}$.
By the measurability and compact-valuedness of $\Gamma(\theta,\cdot)$ in Proposition \ref{prop:subdifferentiable-V} and
the upper semicontinuity and compactness of $M(\theta,\cdot,\cdot)$ in part (i), $\mathcal U(\cdot)$ is measurable and compact-valued.
Then $R$ is a Carath\'eodory mapping, i.e., $R(\cdot,(z,\eta))$ is measurable for each $(z,\eta)$ and $R(x,\cdot)$ is continuous for each $x$. 
By Filippov's measurable selection theorem \cite[Theorem 8.2.10]{aubin1992set}, there exists a measurable selection $(z(\cdot),\eta(\cdot))$ such that
$(z(x),\eta(x))\in \mathcal U(x)$
and
\begin{eqnarray*}
   q(x) &=& R(x,(z(x),\eta(x)))\\
   &=&\nabla_\theta \mathbb{E}_P\left[L(z(x),x,Y,\theta)\mid X=x\right] +\nabla^2_{z\theta}C(z(x),\theta,x)^{\top}\eta(x),\;\forall \inmat{ a.e. } x\in \mathcal{X}.
\end{eqnarray*}
In particular, $z(x)\in\Gamma(\theta,x)$ and $\eta(x)\in M(\theta,x,z(x))$ for a.e.  $x\in \mathcal{X}$,
which establishes the measurable selections for $z(\cdot)$ and $\eta(\cdot)$.
\hfill $\Box$

\subsection{Characterization of the Mordukhovich coderivative for polyhedral feasible set}

The optimality conditions established in Theorem \ref{thm:first-order-necessary-condition} involves the Mordukhovich coderivative term
$D^*\N_\mathcal{Z}\big(z(x), - \nabla_{z}\mathfrak{C}(z(x),\theta,x)) (\eta(x)\big)$. 
In this subsection, we derive an explicit forms 
%characterizations 
of this term in two common settings: (a)~a polyhedral feasible set of the form $\mathcal{Z} = \left\{z \in \mathbb{R}_+^{d_z}: Az \leq b\right\}$ and (b) the nonnegative orthant $\mathcal{Z} = \mathbb{R}_+^{d_z}$ with $A\in \mathbb{R}^{m\times d_z}$ and $b\in \mathbb{R}^m$.
%Recalling t
By the definition of the Mordukhovich coderivative (see Definition~\ref{def:mordukhovich-coderivative}), 
%we have
\begin{eqnarray*}
\zeta \in D^*\N_\Z(\bar z,- \nabla_{z}\mathfrak{C}(z(x),\theta,x)) (\eta) \Longleftrightarrow  (\zeta, - \eta) \in \N_{gph \N_\Z} (\bar z, - \nabla_{z}\mathfrak{C}(z(x),\theta,x)).
\end{eqnarray*}
Thus, computing the Mordukhovich coderivative reduces to determining the limiting normal cone to $\text{gph}\; \N_\Z$.

We begin by recalling some basic notions.
For a polyhedral convex set $\mathcal{Z}$, 
the {\em tangent cone} to $\mathcal{Z}$ at $z$ is defined by 
$$
T_\mathcal{Z}(z):=\left\{ d\in \mathbb{R}^{d_z}: \exists \{z_i\}\subset \mathcal{Z}, z_i\rightarrow z, t_i\downarrow 0, \text{ such that } d=\lim_{i\rightarrow\infty} \frac{z_i-z}{t_i}  \right\},
$$
which is a polyhedral convex set. 
% Recall the normal cone given in Definition \ref{def:normal-cone}.
The {\em critical cone} of $\mathcal{Z}$ at $z$ for a vector $v\in \mathcal{N}_\mathcal{Z}(z)$ is defined by
$$
K_\mathcal{Z}(z,v) :=\left\{d\in T_\mathcal{Z}(z): d\perp v \right\},
$$
which is also a polyhedral convex set.
Let $K^*$ denote the {\em polar cone} of a polyhedral convex cone $K$, i.e.,
$K^*:=\{w\in \mathbb{R}^{d_z}: \langle w,d \rangle\leq 0,\forall d\in K\}$.
A {\em closed face}, denoted by $F$, of a polyhedral convex cone $K$
is %a polyhedral convex cone $K$ is 
defined 
by
%of the form
$$
F=\left\{ d\in K: d\perp w \right\} \text{ for some } w\in K^*.
$$

\begin{lemma} \label{lemma:normal_to_gph_normal_polyhedron}
    Let $\mathcal{Z}\subset \mathbb{R}^{d_z}$ be a nonempty polyhedral convex set, and $(\bar{z},\bar{v})\in \inmat{gph} \mathcal{N}_{\mathcal{Z}}$.
    Denote by $\mathcal{F}(K_\mathcal{Z}(\bar{z},\bar{v}))$ the family of closed faces of the critical cone $K_\mathcal{Z}(\bar{z},\bar{v})$.
    Then the normal cone to the graph of $\mathcal{N}_\mathcal{Z}$ at $(\bar{z},\bar{v})$ takes the form 
    \begin{eqnarray*} 
    % \label{eqn:normal_to_gph_normal_polyhedron}
        \mathcal{N}_{\text{gph}\mathcal{N}_\mathcal{Z}}(\bar{z},\bar{v}) = \bigcup_{\substack{F_1,F_2 \in \mathcal{F}(K_\mathcal{Z}(\bar{z},\bar{v})) \\ F_2\subset F_1}} (F_1-F_2)^*\times (F_1-F_2).
    \end{eqnarray*}
    % where $F_1 - F_2 := \{ d_1 - d_2 \mid d_1 \in F_1,\; d_2 \in F_2\}.$
\end{lemma}

\noindent
\textbf{Proof.}
This result is provided in \cite[Proof of Theorem 2]{dontchev1996characterizations} and \cite[Proof of Proposition 4.4]{poliquin1998tilt}.
For completeness and to facilitate reading, we include a proof in the Appendix \ref{sec:appendix}.
\hfill $\Box$

Based on Lemma \ref{lemma:normal_to_gph_normal_polyhedron},
we provide a explicit characterization of $\mathcal{N}_{\text{gph}\mathcal{N}_\mathcal{Z}}(\bar{z},\bar{v})$ when $\mathcal{Z}=\{z\in \mathbb{R}^{d_z}:Az\leq b\}$ with $a_i\in \mathbb{R}^{d_z}$ for $i=1,\dots, m$, $A = (a_1,\dots,a_m)^{\top}\in \mathbb{R}^{m \times d_z}$ and $b\in \mathbb{R}^m$. 
Let $I(z)$ be the active set of constraints at $z\in \mathcal{Z}$, i.e. $I(z)=\{ i\in \{1\dots,m\}: a_i^{\top}z = b_i \}$.

\begin{proposition} \label{prop:normal_to_azb}
    For $\mathcal{Z}=\{z\in \mathbb{R}^{d_z}:Az\leq b\}$, 
    % $\mathcal{N}_{\text{gph}\mathcal{N}_\mathcal{Z}}(\bar{z},\bar{v})$ is equivalent to 
    \begin{eqnarray}\label{eqn:normal_to_azb}
        \mathcal{N}_{\text{gph}\mathcal{N}_\mathcal{Z}}(\bar{z},-\bar{v})&=&\bigg\{(\zeta,-\eta)\in \mathbb{R}^{d_z}\times \mathbb{R}^{d_z}: \exists \bar{\lambda}\in \mathbb{R}_+^{m} \text{ and } J_1\subset J_2 \subset \{1,\dots,m\}, \nonumber\\
    & &\quad -\bar{v} = A^{\top}\bar{\lambda},\quad \bar{\lambda}_i(a_i^{\top}\bar{z}-b_i) = 0, \quad \forall i=1,\dots,m,\nonumber\\
    & &\quad a_i^{\top} \bar{z} = b_i,\; \bar{\lambda}_i=0,\; \forall i\in J_2,\nonumber\\
    & &\quad \zeta = A^{\top}_{J_2 \backslash J_1}\mu + A^{\top}_{\{i: a_i^{\top}\bar{z} = b_i,\bar{\lambda}_i>0\}\cup J_1}\nu, \mu\in \mathbb{R}_+^{|J_2|-|J_1|}, \nu\in \mathbb{R}^{|J_1|+|\{i: a_i^{\top}\bar{z} = b_i,\bar{\lambda}_i>0\}|},\nonumber\\
    & &\quad a_i^{\top}\eta = 0, \forall i\in \{i: a_i^{\top}\bar{z} = b_i,\bar{\lambda}_i>0\}\cup J_1,\;  a_i^{\top}\eta\geq 0, \forall i\in J_2 \backslash J_1 \bigg\}.
    \end{eqnarray}
\end{proposition}

\noindent
\textbf{Proof.}
By \cite[Proposition 4.1.5]{facchinei2003finite}, the tangent cone and normal cone of $\mathcal{Z}$ at $\bar{z}\in \mathcal{Z}$ are given by
\begin{eqnarray*}
    T_\mathcal{Z}(\bar{z}) &=& \left\{ d\in \mathbb{R}^{d_z}: a_i^{\top} d \leq 0, \forall i\in I(\bar{z}) \right\} \\
    \mathcal{N}_\mathcal{Z}(\bar{z}) &=& \left\{ \sum_{i\in I(\bar{z})} \lambda_ia_i: \lambda \in \mathbb{R}_+^m \right\} = \left\{ A^{\top}\lambda: \lambda \in \mathbb{R}_+^m, \lambda_i(a_i^{\top}\bar{z}-b_i) = 0, i=1,\dots,m \right\}.
\end{eqnarray*}
For $(\bar{z},-\bar{v})\in \text{gph} \mathcal{N}_\mathcal{Z}$, there exists $\bar{\lambda}\geq 0$ such that
\begin{eqnarray*}
    -\bar{v} = A^{\top}\bar{\lambda},\quad \lambda_i(a_i^{\top}\bar{z}-b_i) = 0, \quad \forall i=1,\dots,m.
\end{eqnarray*}
Since $K_\mathcal{Z}(\bar{z},-\bar{v}) = T_\mathcal{Z}(\bar{z}) \cap [-\bar{v}]^\perp$, we have 
\begin{eqnarray*}
    K_\mathcal{Z}(\bar{z},\bar{v}) = \left\{ d\in \mathbb{R}^{d_z} : a_i^{\top} d \leq 0, \forall i\in I(\bar{z}), \langle\bar{v},d\rangle = 0\right\}.
\end{eqnarray*}
Particularly, based on $-\bar{v}= A^{\top}\bar{\lambda}$, we separate $I(\bar{z})$ as 
\begin{eqnarray*}
    I_+(\bar{z}):= \{ i\in I(\bar{z}): \bar{\lambda}_i>0 \}, \quad I_0(\bar{z}):= \{ i\in I(\bar{z}): \bar{\lambda}_i=0 \}.
\end{eqnarray*}
Since 
\begin{eqnarray*}
    \langle\bar{v},d\rangle = -\bar{\lambda}^{\top} A d = -\sum_{i\in I(\bar{z})} \bar{\lambda}_i a_i^{\top} d,
\end{eqnarray*}
we have 
\begin{eqnarray*}
    K_\mathcal{Z}(\bar{z},-\bar{v}) &=&  \left\{ d\in\mathbb{R}^{d_z}:\; a_i^{\top} d \leq 0, \forall i\in I(\bar{z});\; \sum_{i\in I(\bar{z})} \bar{\lambda}_i a_i^{\top} d=0, \bar{\lambda}\geq 0 \right\}\\
    &=& \left\{ d\in\mathbb{R}^{d_z}:\; a_i^{\top} d=0, \forall i\in I_+(\bar{z});\; a_i^{\top} d \leq 0, \forall i\in I_0(\bar{z}) \right\}.
\end{eqnarray*}
Then any closed face of $K_\mathcal{Z}(\bar{z},-\bar{v})$ can be written with $J\in I_0$ as
\begin{eqnarray*}
    F_J:=\left\{ d\in\mathbb{R}^{d_z}:\; a_i^{\top} d=0, \forall i\in I_+(\bar{z})\cup J;\; a_i^{\top} d \leq 0, \forall i\in I_0(\bar{z})\backslash J\right\}.
\end{eqnarray*}
Note that when $J_1\subset J_2\subset I_0(\bar{z})$, $F_{J_2}\subset F_{J_1}\subset K_\mathcal{Z}(\bar{z},\bar{v})$. 
Then 
\begin{eqnarray*}
    F_{J_1} - F_{J_2} 
    =  \left\{  d:\; a_i^{\top}d = 0, \forall i\in I_+(\bar{z})\cup J_1;\;  a_i^{\top}d\leq 0, \forall i\in J_2 \backslash J_1  \right\},
\end{eqnarray*}
and its polar cone is 
\begin{eqnarray*}
    (F_{J_1} - F_{J_2})^*=\left\{ A^{\top}_{J_2 \backslash J_1}\mu + A^{\top}_{I_+(\bar{z})\cup J_1}\nu: \mu\in \mathbb{R}_+^{|J_2|-|J_1|}, \nu\in \mathbb{R}^{|J_1|+|I_+(\bar{z})|} \right\},
\end{eqnarray*}
where $A_{J}$ denotes the submatrix of $A$ with row indices in $J$.
Therefore, by Lemma \ref{lemma:normal_to_gph_normal_polyhedron},
\begin{eqnarray*}
    \mathcal{N}_{\text{gph}\mathcal{N}_\mathcal{Z}}(\bar{z},-\bar{v}) &=& \bigcup_{J_1\subset J_2 \subset I_0(\bar{z}) } (F_{J_1}-F_{J_2})^*\times (F_{J_1}-F_{J_2}) \\
    &=&\bigg\{(\zeta,-\eta)\in \mathbb{R}^{d_z}\times \mathbb{R}^{d_z}: \exists J_1\subset J_2 \subset I_0(\bar{z}), \\
    & &\quad \zeta = A^{\top}_{J_2 \backslash J_1}\mu + A^{\top}_{I_+(\bar{z})\cup J_1}\nu, \mu\in \mathbb{R}_+^{|J_2|-|J_1|}, \nu\in \mathbb{R}^{|J_1|+|I_+(\bar{z})|},\\
    & &\quad a_i^{\top}\eta = 0, \forall i\in I_+(\bar{z})\cup J_1,\;  a_i^{\top}\geq 0, \forall i\in J_2 \backslash J_1 \bigg\}\\
    &=& \bigg\{(\zeta,-\eta)\in \mathbb{R}^{d_z}\times \mathbb{R}^{d_z}: \exists \bar{\lambda}\in \mathbb{R}_+^{m} \text{ and } J_1\subset J_2 \subset \{1,\dots,m\}, \\
    & &\quad -\bar{v} = A^{\top}\bar{\lambda},\quad \bar{\lambda}_i(a_i^{\top}\bar{z}-b_i) = 0, \quad \forall i=1,\dots,m,\\
    & &\quad a_i^{\top} \bar{z} = b_i,\; \bar{\lambda}_i=0,\; \forall i\in J_2,\\
    & &\quad \zeta = A^{\top}_{J_2 \backslash J_1}\mu + A^{\top}_{\{i: a_i^{\top}\bar{z} = b_i,\bar{\lambda}_i>0\}\cup J_1}\nu, \mu\in \mathbb{R}_+^{|J_2|-|J_1|}, \nu\in \mathbb{R}^{|J_1|+|\{i: a_i^{\top}\bar{z} = b_i,\bar{\lambda}_i>0\}|},\\
    & &\quad a_i^{\top}\eta = 0, \forall i\in \{i: a_i^{\top}\bar{z} = b_i,\bar{\lambda}_i>0\}\cup J_1,\;  a_i^{\top}\eta\geq 0, \forall i\in J_2 \backslash J_1 \bigg\}
\end{eqnarray*}
which is exactly \eqref{eqn:normal_to_azb}.
\hfill $\Box$

Combining Theorem \ref{thm:first-order-necessary-condition} and Proposition \ref{prop:normal_to_azb}, 
we are ready to 
%we provide 
present the first-order 
necessary optimality condition of the ILO model \eqref{eq:ILO-3} when $\mathcal{Z}$ has a polyhedral structure. 
%with the polyhedral feasible set $\mathcal{Z}$.

\begin{theorem}
%Consider the case where 
Let $\mathcal{Z}=\{z\in \mathbb{R}^{d_z}:Az\leq b\}$ and the assumptions in Theorem \ref{thm:first-order-necessary-condition} hold. 
If $(\theta,z(\cdot))$ is a local optimal solution of problem \eqref{eq:ILO-3}, then
% Then the first-order necessary optimality condition of the ILO model \eqref{eq:ILO-3} is that 
there exist $\eta:\mathcal{X}\rightarrow\mathbb{R}^{d_z}$, $\zeta:\mathcal{X}\rightarrow\mathbb{R}^{d_z}$, $\lambda: \mathcal{X} \rightarrow \mathbb{R}_+^{m}$, $J_1: \mathcal{X} \rightarrow 2^{\{1,\dots,m\}}$, and $J_2: \mathcal{X} \rightarrow 2^{\{1,\dots,m\}}$, such that
\begin{eqnarray*}
        \begin{cases}
            0 \in  \mathbb{E}_P \left[ 
            \nabla_\theta \mathbb{E}_P\left[L(z(X),X,Y,\theta)\mid X\right] + 
            \nabla_{z\theta}^2 \mathfrak{C}(z(X),\theta,X)^{\top} \eta(X) \right] + \mathcal{N}_\Theta(\theta),\\
            0 \in \nabla_{z} \mathbb{E}_P\left[L(z(X),X,Y,\theta)\mid X=x\right] + \nabla_{zz}^2 \mathfrak{C}(z(x),\theta,x)^{\top} \eta(x)  + \zeta(x), \; \forall x\in \mathcal{X},\\
            - \nabla_{z}\mathfrak{C}(z(x),\theta,x) = A^{\top}\lambda(x),\quad \lambda_i(x)(a_i^{\top}z(x)-b_i) = 0, \quad \forall i=1,\dots,m,\;  x\in \mathcal{X},\\
            J_1(x)\subset J_2(x) \subset I_0(z(x)), \; \forall x\in \mathcal{X}, \\
            \zeta(x) = A^{\top}_{J_2(x) \backslash J_1(x)}\mu(x) + A^{\top}_{I_+(z(x))\cup J_1(x)}\nu(x),\; \forall x\in \mathcal{X},\\ 
            \mu(x)\in \mathbb{R}_+^{|J_2(x)|-|J_1(x)|}, \nu(x)\in \mathbb{R}^{|J_1(x)|+|I_+(z(x))|},\; \forall x\in \mathcal{X},\\
            a_i^{\top}\eta(x) = 0, \forall i\in I_+(z(x))\cup J_1(x),\;  a_i^{\top}\eta(x)\geq 0, \forall i\in J_2(x) \backslash J_1(x),\;  x\in \mathcal{X},\\
            0 \in \nabla_z \mathfrak{C}(z(x),\theta,x) +{\cal N}_{\cal Z}(z), \; \forall x\in \mathcal{X},
        \end{cases}
\end{eqnarray*}
% where $L(x),I_+(x),I_0(x)$ are defined in \eqref{e.qn:label-first-orthant}.
where $I_+(z(x)):= \{ i\in \{1,\dots,m\}: \lambda_i(x)>0, a_i^{\top}z(x) = b_i \}$, and $I_0(z(x)):= \{ i\in \{1,\dots,m\}: \lambda_i(x)=0, a_i^{\top}z(x) = b_i \}$,
or equivalently,
\begin{eqnarray*}
        \begin{cases}
            0 \in  \mathbb{E}_P \left[ 
            \nabla_\theta \mathbb{E}_P\left[L(z(X),X,Y,\theta)\mid X\right] + 
            \nabla_{z\theta}^2 \mathfrak{C}(z(X),\theta,X)^{\top} \eta(X) \right] + \mathcal{N}_\Theta(\theta),\\
            0 \in \nabla_{z} \mathbb{E}_P\left[L(z(X),X,Y,\theta)\mid X=x\right] + \nabla_{zz}^2 \mathfrak{C}(z(x),\theta,x)^{\top} \eta(x)  + \zeta(x),\; \forall x\in \mathcal{X},\\
            - \nabla_{z}\mathfrak{C}(z(x),\theta,x) = A^{\top}\lambda(x),\quad \lambda_i(x)(a_i^{\top}z(x)-b_i) = 0, \quad \forall i=1,\dots,m,\; x\in \mathcal{X},\\
            J_1(x)\subset J_2(x),\; \forall x\in \mathcal{X},\\
            a_i^{\top} z(x) = b_i,\; \lambda_i(x)=0,\; \forall i\in J_2(x),\; \forall x\in \mathcal{X},\\
            \zeta(x) = A^{\top}_{J_2(x) \backslash J_1(x)}\mu(x) + A^{\top}_{\{i: a_i^{\top}z(x) = b_i,\lambda_i(x)>0\}\cup J_1(x)}\nu(x), \; \forall x\in \mathcal{X},\\
            \mu(x)\in \mathbb{R}_+^{|J_2(x)|-|J_1(x)|}, \nu(x)\in \mathbb{R}^{|J_1(x)|+|\{i: a_i^{\top}z(x) = b_i,\lambda_i(x)>0\}|},\; \forall x\in \mathcal{X},\\
            a_i^{\top}\eta(x) = 0, \forall i\in \{i: a_i^{\top}z(x) = b_i,\lambda_i(x)>0\}\cup J_1(x),\;  a_i^{\top}\eta(x)\geq 0, \forall i\in J_2(x) \backslash J_1(x),\;  x\in \mathcal{X},\\
            0 \in \nabla_z \mathfrak{C}(z(x),\theta,x) +{\cal N}_{\cal Z}(z),\; \forall x\in \mathcal{X}.
        \end{cases}
\end{eqnarray*}
% where $L(x),I_+(x),I_0(x)$ are defined in \eqref{e.qn:label-first-orthant}.
\end{theorem}

% {\color{red}HX: 1. add remark; 2. Examples, 3. Meaurability, integrability?? }
% {\color{blue} Yuan: the measurability is discussed in aforementioned blue part.}

{\color{black}

Next, we consider a particularly interesting case relevant to our discussion when $\mathcal{Z}=\mathbb{R}_+^{d_z}$.
By setting $A=-I^{d_z\times d_z}$ and $b=0$ in Proposition \ref{prop:normal_to_azb}, we have the following result.

\begin{corollary} \label{prop:coderivative-first-orthant}
For any $(z, - v) \in gph N_{\mathbb{R}_+^{d_z}}$ and $x\in \mathcal{X}$, let
\begin{subequations} \label{eqn:label-first-orthant}
    \begin{eqnarray}
L
&:=&\{ i\in \{1,2,\dots, d_z\}: z_i >0, v_i =0\},\\
I_+
&:=&\{ i\in \{1,2,\dots, d_z\}: z_i =0, v_i >0\},\\
I_0
&:=&\{ i\in \{1,2,\dots, d_z\}: z_i =0, v_i =0\}.
\end{eqnarray}
\end{subequations}
Then 
\begin{eqnarray}\label{eqn:coderivative-first-orthant}
\N_{gph \N_{\mathbb{R}_+^{d_z}}}(z,-v) &=&\{(\zeta,-\eta) \in\mathbb{R}^{d_z}\times \mathbb{R}^{d_z}: \zeta_{L}=0, \eta_{I_+}=0,\nonumber\\
&& \forall i \in I_0, \mbox{ either } \zeta_i<0, \eta_i<0 \mbox{ or } \zeta_i \eta_i=0\}.
\end{eqnarray}
\end{corollary}

\noindent
\textbf{Proof.}
In Proposition \ref{prop:normal_to_azb}, 
we consider the case with $\mathcal{Z}=\mathbb{R}_+^{d_z}$ by setting $A=-I^{d_z\times d_z}$ and $b=0$.
Then 
\begin{eqnarray*}
        \mathcal{N}_{\text{gph}\mathcal{N}_\mathcal{Z}}(\bar{z},-\bar{v})&=&\bigg\{(\zeta,-\eta)\in \mathbb{R}^{d_z}\times \mathbb{R}^{d_z}: \exists \bar{\lambda}\in \mathbb{R}_+^{m} \text{ and } J_1\subset J_2 \subset \{1,\dots,m\}, \nonumber\\
    & &\quad \bar{v} = \bar{\lambda},\quad \bar{\lambda}_i \bar{z}_i= 0, \quad \forall i=1,\dots,m,\nonumber\\
    & &\quad \bar{z}_i = 0,\; \bar{\lambda}_i=0,\; \forall i\in J_2,\nonumber\\
    & &\quad \zeta = \big(-I^{d_z\times d_z}\big)^{\top}_{J_2 \backslash J_1}\mu + \big(-I^{d_z\times d_z}\big)^{\top}_{\{i: \bar{z}_i = 0,\bar{\lambda}_i>0\}\cup J_1}\nu,\\
    && \quad\mu\in \mathbb{R}_+^{|J_2|-|J_1|}, \nu\in \mathbb{R}^{|J_1|+|\{i: \bar{z}_i = 0, \bar{\lambda}_i>0\}|},\nonumber\\
    & &\quad \eta_i = 0, \forall i\in \{i: \bar{z}_i = 0,\bar{\lambda}_i>0\}\cup J_1,\;  \eta_i\leq 0, \forall i\in J_2 \backslash J_1 \bigg\}.
    \end{eqnarray*}
Eliminating $\lambda$,
\begin{eqnarray*}
        \mathcal{N}_{\text{gph}\mathcal{N}_\mathcal{Z}}(\bar{z},-\bar{v})&=&\bigg\{(\zeta,-\eta)\in \mathbb{R}^{d_z}\times \mathbb{R}^{d_z}: \exists  J_1\subset J_2 \subset \{1,\dots,m\}, \nonumber\\
    & &\quad \bar{z}_i = 0,\; \bar{v}_i=0,\; \forall i\in J_2,\nonumber\\
    & &\quad \zeta = \big(-I^{d_z\times d_z}\big)^{\top}_{J_2 \backslash J_1}\mu + \big(-I^{d_z\times d_z}\big)^{\top}_{\{i: \bar{z}_i = 0,\bar{v}_i>0\}\cup J_1}\nu, \\
    & & \quad \mu\in \mathbb{R}_+^{|J_2|-|J_1|}, \nu\in \mathbb{R}^{|J_1|+|\{i: \bar{z}_i = 0, \bar{\lambda}_i>0\}|},\nonumber\\
    & &\quad \eta_i = 0, \forall i\in \{i: \bar{z}_i = 0,\bar{v}_i>0\}\cup J_1,\;  \eta_i\leq 0, \forall i\in J_2 \backslash J_1 \bigg\}.
    \end{eqnarray*}
For any $i\in L:=\{ i\in \{1,2,\dots, d_z\}: z_i >0, v_i =0\}$, we have $i\notin \{i: \bar{z}_i = 0,\bar{v}_i>0\}\cup J_2$, and then $\zeta_i=0$.
For $i\in I_+ := \{ i\in \{1,2,\dots, d_z\}: z_i =0, v_i >0\}$, we have $\zeta_i=-\nu_i\in \mathbb{R}$ (which means there is no constraint on $\zeta_i$) and $\eta_i=0$.
Finally, we examine the last case with $i\in I_0 := \{ i\in \{1,2,\dots, d_z\}: z_i =0, v_i =0\}$, where
if $i\notin J_2$, we require $\zeta_i=0$; if $i\in J_2\backslash J_1$, we have $\zeta_i=-\mu_i\leq 0$ and $\eta_i \leq 0$; if $i\in J_1$, we have $\zeta_i=-\nu_i\in \mathbb{R}$ (no constraint on $\zeta_i$) and $\eta_i=0$.
We can summarize these three sub-cases and rewrite it as either $\zeta_i\eta_i=0$ or $\zeta_i< 0,\eta_i<0$.
Then we arrive at \eqref{eqn:coderivative-first-orthant}.
% Therefore, Proposition \ref{prop:coderivative-first-orthant} can be considered a corollary of Proposition \ref{prop:normal_to_azb}.
\hfill $\Box$

In this case, the necessary optimality condition can be expressed in the following form.

\begin{corollary} \label{corollary:opt-cond-first-orthant}
    Consider the case where $\mathcal{Z}:= \mathbb{R}_+^{d_z}$.
Let the assumptions in Theorem \ref{thm:first-order-necessary-condition} hold. 
If $(\theta,z(\cdot))$ is a local optimal solution of problem \eqref{eq:ILO-3}, 
% Then the first-order necessary optimality conditions of the ILO model \eqref{eq:ILO-3} is that 
then there exist $\eta:\mathcal{X}\rightarrow\mathbb{R}^{d_z}$ and $\zeta:\mathcal{X}\rightarrow\mathbb{R}^{d_z}$, such that
\begin{eqnarray*}
        \begin{cases}
            0 \in  \mathbb{E}_P \left[ 
            \nabla_\theta \mathbb{E}_P\left[L(z(X),X,Y,\theta)\mid X\right] + 
            \nabla_{z\theta}^2 \mathfrak{C}(z(X),\theta,X)^{\top} \eta(x) \right] + \mathcal{N}_\Theta(\theta),\\
            0 \in \nabla_{z} \mathbb{E}_P\left[L(z(X),X,Y,\theta)\mid X=x\right] + \nabla_{zz}^2 \mathfrak{C}(z(x),\theta,x)^{\top} \eta(x)  + \zeta(x), \forall x\in \mathcal{X},\\
            \zeta_{L(x)}(x) = 0, \; \eta_{I_+(x)}(x) = 0, \;
            \forall i \in I_0(x), \mbox{ either } \zeta_i(x)<0, \eta_i(x)<0 \mbox{ or } \zeta_i(x) \eta_i(x)=0, \forall x\in \mathcal{X},\\
            0 \in \nabla_z \mathfrak{C}(z(x),\theta,x) +{\cal N}_{\cal Z}(z(x)), \forall x\in \mathcal{X},
        \end{cases}
\end{eqnarray*}
where $L(x),I_+(x),I_0(x)$ are defined in \eqref{eqn:label-first-orthant} by letting $v= \nabla_z \mathfrak{C}(z(x),\theta,x)$ for each $x\in \mathcal{X}$.
\end{corollary}

}

The first-order optimality conditions can be extended to the case when $\mathcal{Z}$ is a simplex in %the first orthant of 
$\R^{d_z}$, the next corollary states this.

\begin{corollary} \label{corollary:simplex}
    For $z\in \mathcal{Z}=\{z\in \mathbb{R}^{d_z}: z\geq 0, 1^\top z\leq 1\}$, let $L:= \{i=\{1,\dots,d_z\}: z_i>0\}$.
    If $1^\top z< 1$, then $I_+:=\{i=\{1,\dots,d_z\}: z_i=0, v_i>0\}$, $ I_0:=\{i=\{1,\dots,d_z\}: z_i=0, v_i = 0\}$;
    otherwise, if $1^\top z=1$, then $L$ is nonempty, $v_j=v_k$ for all $j,k\in L$, and let $\tau:=-v_i$ for any $i\in L$, $I_+:=\{i=\{1,\dots,d_z\}: z_i=0, v_i>-\tau\}$, $ I_0:=\{i=\{1,\dots,d_z\}: z_i=0, v_i = -\tau\}$.
    Then
    \begin{eqnarray*} 
    % \label{eqn:simplex-coderivative}
        \mathcal{N}_{\text{gph}\mathcal{N}_\mathcal{Z}}(z,-v)&=&\bigg\{(\zeta,-\eta)\in \mathbb{R}^{d_z}\times \mathbb{R}^{d_z}: \exists \beta\in \mathbb{R},\nonumber\\
        && \quad \beta (1-\boldsymbol{1}^\top z) = 0, \; \tau(\boldsymbol{1}^\top \eta) = 0, \nonumber\\
        && \quad \text{either } \beta>0, \boldsymbol{1}^\top \eta >0, \text{ or } \beta \boldsymbol{1}^\top \eta = 0,\nonumber \\
        &&\quad \zeta_L = \beta, \; \eta_{I_+}=0,\nonumber\\
        &&\quad \forall i\in I_0, \text{ either } \zeta_i<\beta, \eta_i<0, \text{ or } (\zeta_i - \beta)\eta_i = 0 \bigg\}.
    \end{eqnarray*}
\end{corollary}

\noindent
\textbf{Proof.}
This corollary can be obtained by considering the polyhedron feasible set in Proposition \ref{prop:normal_to_azb} by setting
%\begin{eqnarray*}
 $   A= 
    \begin{bmatrix}
        -I^{d_z\times d_z}\\
        \boldsymbol{1}^\top
    \end{bmatrix}, %\qquad
    b=
    \begin{bmatrix}
        \boldsymbol{0}\\
        1
    \end{bmatrix},
$
%\end{eqnarray*}
where $\boldsymbol{0},\boldsymbol{1}$ are $d_z$-dimensional vectors.
See Appendix \ref{sec:proof-corollary-5.1} for details.
\hfill $\Box$

\section{First-order Necessary Optimality Condition: Nonconvex Case}
% lower-level problem \eqref{eqn:lower-level-problem-A} is nonconvex}
%\section{Necessary optimality conditions}
%: combined SMPEC and the value function approach}
\label{eqn:nonconvex}

In this section, we extend the discussions in the preceding section to 
the case where $\mathfrak{C}(z,\theta, x)$ is nonconvex in $z$.
In that case, problem \eqref{eq:ILO-2} 
%may no longer be
is no longer equivalent to problem \eqref{eq:ILO-3} in general. 
Some remedies are needed for the latter before 
we derive the first-order optimality conditions for the former.
Similar to deterministic bilevel program case, we adopt
%% In this case, the KKT approach is generally not reliable, because many stationary points and local maximizers will also satisfy these conditions.
%In this section, we follow the idea in 
the calmness conditions introduced by
Ye and Zhu \cite{ye2010new}.
%and introduce the value function of the decision-making problem to characterize the optimality condition.
Throughout this section, we make the following assumption.

\begin{assumption} \label{assumption:convex-smooth2}
    For any $(\theta,x)\in \Theta\times \mathcal{X}$, $\mathfrak{C}$ is continuously differentiable in $z$. 
    The feasible set $\mathcal{Z}$ is a closed, convex subset of $\mathbb{R}^{d_z}$.
    % The feasible set $\mathcal{Z}$ is a compact, convex subset of $\mathbb{R}^{d_z}$.
\end{assumption}

\subsection{
%Equivalent problems and u
Uniform/stochastic partial calmness.}

To ease the exposition, we recall that the problem \eqref{eq:ILO-1} is formulated as follows: 
% equivalent to
\begin{subequations}
\label{eq:ILO-2-rewrite}
    \begin{eqnarray}
    \min_{\theta\in\Theta,z(\cdot)\in \mathfrak{M}}
        && 
        % \mathbb{E}_{P_X}\left[\mathbb{E}_{P_{Y|X}}[L(z^*(\theta,x),X,Y)]\right] = 
        \mathbb{E}_{P} \left[
        %\min_{z^*(\theta,x)}
        L(z(X),X,Y,\theta)\right] \label{eqn:ILO-2-rewrite-upper-level-problem}\\
        \inmat{s.t.} && z(x)\in \arg\min_{z\in \mathcal{Z}} \mathfrak{C}(z,\theta,x),\; \inmat{for a.e.} \;x\in \mathcal{X}, \label{eq:ILO-2-rewrite-lower-level}
\end{eqnarray}
\end{subequations}
where $\mathfrak{M}$ denotes the set of measurable functions $z:\mathcal{X}\rightarrow \mathcal{Z}$. 
Let $\mathfrak{C}^*(\theta,x) = \min_{z\in \mathcal{Z}}  \mathfrak{C}(z,\theta,x)$ and $\mathcal{Z}^*(\theta,x)=\arg\min_{z\in \mathcal{Z}} \mathfrak{C}(z,\theta,x)$ denote the optimal value and the set of optimal solutions for the lower-level problem \eqref{eq:ILO-2-rewrite-lower-level} for $x\in \mathcal{X},\theta\in \Theta$.
%The combined SMPEC and the value function approach consider the following combined problem:
The two-stage bilevel program above can be equivalently written as
\begin{subequations}
\label{eqn:CP}
    \begin{eqnarray}
        % \text{(CP-1)} \qquad
        \min_{\theta\in\Theta,z(\cdot)\in \mathfrak{M}}
        && 
        % \mathbb{E}_{P_X}\left[\mathbb{E}_{P_{Y|X}}[L(z^*(\theta,x),X,Y)]\right] = 
        \mathbb{E}_{P} \left[
        %\min_{z^*(\theta,x)}
        L(z(X),X,Y,\theta)\right]\\
        \inmat{s.t.} && 0\in \nabla_z \mathfrak{C}(z(x),\theta,x) +{\cal N}_{\cal Z}(z(x)),\; \inmat{for a.e.} \;x\in \mathcal{X}, \label{eqn:CP-ec} \\
        && \mathfrak{C}(z(x),\theta,x) - \mathfrak{C}^*(\theta,x)\leq 0, \; \inmat{for a.e.} \;x\in \mathcal{X}, \label{eqn:CP-value}
    \end{eqnarray}
\end{subequations}
and further as
\begin{subequations}
\label{eqn:CP-2}
    \begin{eqnarray}
        % \text{(CP-2)} \qquad
        \min_{\theta\in\Theta,z(\cdot)\in \mathfrak{M}}
        && 
        % \mathbb{E}_{P_X}\left[\mathbb{E}_{P_{Y|X}}[L(z^*(\theta,x),X,Y)]\right] = 
        \mathbb{E}_{P} \left[
        %\min_{z^*(\theta,x)}
        L(z(X),X,Y,\theta)\right]\\
        \inmat{s.t.} && 0\in \nabla_z \mathfrak{C}(z(x),\theta,x) +{\cal N}_{\cal Z}(z(x)),\; \inmat{for a.e.} \;x\in \mathcal{X}, \label{eqn:CP2-ec} \\
        && \mathbb{E}\left[\mathfrak{C}(z(x),\theta,x) - \mathfrak{C}^*(\theta,x) \right] \leq 0. \label{eqn:CP2-value}
    \end{eqnarray}
\end{subequations}
While the value function constraint \eqref{eqn:CP-value}
effectively avoids potential discrepancy between 
constraint \eqref{eq:ILO-2-rewrite-lower-level} and 
constraint \eqref{eqn:CP-ec}, 
it also leads to abnormal multipliers and the failure of any constraint qualification. 
To tackle the issue, we will need to introduce 
the so-called partial calmness condition. 

\begin{definition}[Uniform/stochastic partial calmness]
    Let $(\bar{\theta},\bar{z}(\cdot))$ be a local optimal solution of problem \eqref{eq:ILO-2-rewrite}; thus, it is also a local solution of \eqref{eqn:CP} and \eqref{eqn:CP-2}.  
    We say that problem \eqref{eqn:CP-2} is uniformly partially calm at $(\bar{\theta},\bar{z}(\cdot))$ if there exists $\mu> 0$ such that $(\bar{\theta},\bar{z}(\cdot))$ is a local optimal solution of the following partially penalized problem:
    \begin{subequations}
\label{eqn:CP-mu}
    \begin{eqnarray}
        % \text{(CP)}_\mu \qquad
        \min_{\theta\in\Theta,z(\cdot)\in \mathfrak{M}}
        && 
        % \mathbb{E}_{P_X}\left[\mathbb{E}_{P_{Y|X}}[L(z^*(\theta,x),X,Y)]\right] = 
        \mathbb{E}_{P} \left[
        %\min_{z^*(\theta,x)}
        L(z(X),X,Y,\theta)\right] + \mu \mathbb{E}_P\left[ \mathfrak{C}(z(x),\theta,X) - \mathfrak{C}^*(\theta,X) \right] \\
        \inmat{s.t.} && 0\in \nabla_z \mathfrak{C}(z(x),\theta,x) +{\cal N}_{\cal Z}(z(x)),\; \forall \inmat{ a.e. } x\in \mathcal{X}. \label{eqn:cp-mu-ec}
    \end{eqnarray}
\end{subequations}
Likewise, we say that problem \eqref{eqn:CP} is stochastically partially calm at $(\bar{\theta},\bar{z}(\cdot))$ if there exists a measurable function $\mu:\mathcal{X}\rightarrow \mathbb{R}_+$ such that $(\bar{\theta},\bar{z}(\cdot))$ is a locally optimal solution of the following partially penalized problem:
    \begin{subequations}
\label{eqn:CP2-mu}
    \begin{eqnarray}
        % \text{(CP)}_{\mu(\cdot)} \qquad
        \min_{\theta\in\Theta,z(\cdot)\in \mathfrak{M}}
        && 
        % \mathbb{E}_{P_X}\left[\mathbb{E}_{P_{Y|X}}[L(z^*(\theta,x),X,Y)]\right] = 
        \mathbb{E}_{P} \left[
        %\min_{z^*(\theta,x)}
        L(z(X),X,Y,\theta)\right] + \mathbb{E}_P\left[ \mu(X) \left(\mathfrak{C}(z(X),\theta,X) - \mathfrak{C}^*(\theta,X) \right) \right]\qquad \\
        \inmat{s.t.} && 0\in \nabla_z \mathfrak{C}(z(x),\theta,x) +{\cal N}_{\cal Z}(z(x)),\; \forall \inmat{ a.e. } x\in \mathcal{X}. \label{eqn:cp-mu-ec2}
    \end{eqnarray}
\end{subequations}
\end{definition}
% {\color{red}HX: Discuss sufficient conditions which ensure the partial calmness}
% {\color{blue}

Partial calmness is often either unattainable or challenging to verify in practical settings. 
Many studies aim to relax the partial calmness assumption and propose alternative conditions that are more easily verifiable.
For instance,
Mehlitz et al. \cite{mehlitz2021note} examine this problem and offer a detailed analysis for the specific case in which the lower-level problem is linear. 
Ye and Zhu \cite{ye2010new} relax the partial calmness condition to more verifiable conditions which are characterized by the directional derivative of the objective function and the linearization cone of the feasible region at the reference point.
Interested readers are referred to \cite{ke2022generic,mehlitz2021note,royset2024stability,ye2010new} for detailed discussions on the verification of calmness conditions.
Since verifiable sufficient conditions for (stochastic/uniform) partial calmness are highly problem-dependent, we adopt this as a standing assumption in this section to maintain focus on our main line of analysis.

The next proposition states that uniform partial calmness implies stochastic partial calmness, while conversely, stochastic partial calmness implies uniform partial calmness under essential boundedness for $\mu(\cdot)$.

\begin{proposition} \label{prop:uniform/stochastic-PC}
\begin{itemize}
    \item[(i)] If \eqref{eqn:CP-2} is uniformly partially calm at $(\bar{\theta},\bar{z}(\cdot))$ with $\bar{\mu}> 0$, then \eqref{eqn:CP} is uniformly partially calm at $(\bar{\theta},\bar{z}(\cdot))$ with $\mu(x)=\bar{\mu}$ for a.e. $x\in \mathcal{X}$.

    \item[(ii)] If \eqref{eqn:CP} is stochastically partially calm at $(\bar{\theta},\bar{z}(\cdot))$ with a measurable function $\bar{\mu}:\mathcal{X}\rightarrow \mathbb{R}_+$, and $\bar{\mu}(\cdot)$ is essentially bounded, i.e., 
    $$\operatorname*{ess\,sup}\bar{\mu}(\cdot)
    :=\inf\Big\{M\in\mathbb R:\ \mu(x)\le M\ \text{for } x\in \mathcal{X },\text{ a.e.}\Big\}< \infty,$$ 
    then \eqref{eqn:CP-2} is uniformly partially calm at $(\bar{\theta},\bar{z}(\cdot))$ with $\mu^{b}:= \operatorname*{ess\,sup}\bar{\mu}(\cdot)$.
\end{itemize}
\end{proposition}

\noindent
\textbf{Proof.}
Part (i) holds directly by the definitions of uniform partial calmness and stochastic partial calmness.
For part (ii), $(\bar{\theta},\bar{z}(\cdot))$ is feasible to \eqref{eqn:CP2-mu} and thus is also feasible to \eqref{eqn:CP-mu}.
For any feasible solution $(\theta, z(\cdot))$ near $(\bar{\theta},\bar{z}(\cdot))$, we have $\mathfrak{C}(z(x),\theta,x) \geq \mathfrak{C}^*(\theta,x)$, and thus 
\begin{eqnarray*}
    &&\mathbb{E}_{P} \left[L(z(X),X,Y,\theta)\right] + \mu^b \mathbb{E}_P\left[ \mathfrak{C}(z(x),\theta,x) - \mathfrak{C}^*(\theta,x) \right] \\
    &=& \mathbb{E}_{P} \left[L(z(X),X,Y,\theta)\right] + \mathbb{E}_P\left[ \bar{\mu}(X) \left( \mathfrak{C}(z(x),\theta,x) - \mathfrak{C}^*(\theta,x) \right) \right]\\
    && + \mathbb{E}_P\left[ (\mu^b - \bar{\mu}(X)) \left( \mathfrak{C}(z(x),\theta,x) - \mathfrak{C}^*(\theta,x) \right) \right]\\
    &\geq & \mathbb{E}_{P} \left[L(z(X),X,Y,\theta)\right] + \mathbb{E}_P\left[ \bar{\mu}(X) \left( \mathfrak{C}(z(x),\theta,x) - \mathfrak{C}^*(\theta,x) \right) \right]\\
    &\geq &\mathbb{E}_{P} \left[L(\bar{z}(X),X,Y,\theta)\right] + \mathbb{E}_P\left[ \bar{\mu}(X) \left( \mathfrak{C}(\bar{z}(x),\bar{\theta},x) - \mathfrak{C}^*(\bar{\theta},x) \right) \right]\\
    &\geq &\mathbb{E}_{P} \left[L(\bar{z}(X),X,Y,\theta)\right] + \mu^b \mathbb{E}_P\left[ \mathfrak{C}(\bar{z}(x),\bar{\theta},x) - \mathfrak{C}^*(\bar{\theta},x)  \right],
\end{eqnarray*}
where the first inequality comes from $\mu^b\geq \bar{\mu}(x)$ for $x\in \mathcal{X}$ a.e., the second inequality follows from the definition of stochastic partial calmness, and the last inequality is obtained by $\mathfrak{C}(\bar{z}(x),\bar{\theta},x) = \mathfrak{C}^*(\bar{\theta},x)$ for $x\in \mathcal{X}$, a.e..
Therefore, $(\bar{\theta},\bar{z}(\cdot))$ is also a local optimal solution for \eqref{eqn:CP-mu} with $\mu^b$, and we immediately arrive at the result.
\hfill $\Box$

\subsection{Necessary optimality conditions under uniform/stochastic partial calmness} \label{sec:opt-cond-spc}

Proposition \ref{prop:uniform/stochastic-PC} indicates that stochastic partial calmness is a weaker condition than uniform partial calmness. Consequently, in this part, we adopt the stochastic partial calmness assumption.
% and concentrate on problem \eqref{eqn:CP2-mu}.

\begin{assumption} \label{assumption:stochastic-pc}
    Problem \eqref{eqn:CP-2} is stochastically partially calm at $(\bar{\theta},\bar{z}(\cdot))$ with a measurable function $\bar{\mu}:\mathcal{X}\rightarrow \mathbb{R}_+$.
\end{assumption}

Under Proposition \ref{prop:uniform/stochastic-PC} and Assumption \ref{assumption:stochastic-pc}, we restrict our attention to problem \eqref{eqn:CP2-mu}.
Let 
\begin{eqnarray*}
    L_{\mu(x)}(z,x,y,\theta):=L(z,x,y,\theta) + \mu(x) \left(\mathfrak{C}(z,\theta,x) - \mathfrak{C}^*(\theta,x) \right).
\end{eqnarray*}
By the interchangeability principle, problem \eqref{eqn:CP2-mu} can be equivalently written as follows:
\begin{subequations}
\label{eqn:CP2-mu-vfunc}
    \begin{eqnarray}
        % \text{(CP)}_{\mu(\cdot)} \qquad
        \min_{\theta\in\Theta}
        && 
        % \mathbb{E}_{P_X}\left[\mathbb{E}_{P_{Y|X}}[L(z^*(\theta,x),X,Y)]\right] = 
        \mathbb{E}_{P} \left[
        v_{\mu(x)}(\theta,x)
        \right]  \\
        \inmat{s.t.} && v_{\mu(x)}(\theta,x) := \min_{z} \mathbb{E}_P[L_{\mu(x)}(z,X,Y,\theta)|X=x]\\
        &&\;\qquad \qquad\qquad \inmat{s.t.} \;0\in \nabla_z \mathfrak{C}(z,\theta,x) +{\cal N}_{\cal Z}(z),\; \forall x\in \mathcal{X}. \label{eqn:cp-mu-vfuncec-2}
    \end{eqnarray}
\end{subequations}
This problem closely resembles problem \eqref{eq:ILO-3}, with the distinction that the second-stage objective function is penalized through a value-function constraint involving the parameter $\mu(x)$.
We employ the method outlined in Section \ref{sec:kkt-approach} to obtain the necessary optimality conditions for problem \eqref{eqn:CP2-mu-vfunc} while adjusting the underlying assumptions and proof technique accordingly.

First, we need to characterize the properties of the optimal value function $\mathfrak{C}^*(\theta,x)$. 
To this end, we require the existence of a uniformly bounded optimal solution.

\begin{assumption} \label{assumption:infcompactness-feasible-set}
    % For any $(\theta,x)\in \Theta\times \mathcal{X}$, 
    There exist a constant $\beta$ and a compact set $\mathcal{Z}_\beta\in \mathcal{Z}$ such that the lower level set $\left\{ z\in\mathcal{Z}: \mathfrak{C}(z,\theta,x)\leq \beta \right\}$ is nonempty and is contained in $\mathcal{Z}_\beta$ for all $(\theta,x)\in \Theta\times \mathcal{X}$.
\end{assumption}

This uniform inf-compactness condition is relaxed from an explicit assumption of the compactness of the feasible set $\mathcal{Z}$.
This kind of condition is widely used in optimization literature to ensure the well-posedness of the optimization problem, see e.g. \cite{guo2021existence}.
The next proposition guarantees continuous differentiability in a neighborhood of the optimal solution and implies that the subdifferential can be represented as a finite convex combination of at most $d_\theta + 1$ elements at a time. 
This representation significantly simplifies the characterization of the necessary optimality conditions.

\begin{proposition} \label{prop:Cstar-lipschitz-differentiable}
    Suppose Assumption \ref{assumption:convex-smooth2} holds. 
    % Assume in addition that $\mathcal{Z}$ is compact.
    % , and $\mathfrak{C}$ is continuously differentiable in $(z,\theta)$.
    Then for $\bar{\theta}\in \Theta$ and $x\in \mathcal{X}$,
    \begin{itemize}
        \item[(i)] $\mathfrak{C}^*(\theta,x)$ is locally Lipschitz around $\bar{\theta}$, its Clarke subdifferential satisfies
    \begin{eqnarray} \label{eqn:clarke-subdiff-value-function}
        \partial_\theta \mathfrak{C}^*(\bar{\theta},x) = \text{co} \left\{ \nabla_\theta \mathfrak{C}(z,\bar{\theta},x): z \in \mathcal{Z}^*(\bar{\theta}, x) \right\},
    \end{eqnarray}
    and its directional derivative is given by
    \begin{eqnarray*}
        D_\theta \mathfrak{C}^*(\bar{\theta},x;d) = \min_{z\in \mathcal{Z}^*(\bar{\theta},x)} \langle \nabla_\theta \mathfrak{C}(z,\bar{\theta},x),d \rangle, \; \forall d\in \mathbb{R}^n;
    \end{eqnarray*}

        \item[(ii)] for any vector $w\in \partial_\theta \mathfrak{C}^*(\bar{\theta},x)\subset \mathbb{R}^{d_\theta}$,
        there exist points $z_1,\dots,z_{d_\theta+1}\in \mathcal{Z}^*(\bar{\theta},x)$ and $\lambda_1,\dots,\lambda_{d_\theta+1}\geq0$ with $\sum_{i=1}^{d_\theta+1} \lambda_i=1$ such that
        \begin{eqnarray*} 
        % \label{eqn:finite-combination-representation}
            w = \sum_{i=1}^{d_\theta+1} \lambda_i \nabla_\theta \mathfrak{C}(z_i,\bar{\theta},x);
        \end{eqnarray*}
    \end{itemize}

    \item[(iii)] for each $x\in \mathcal{X}$, $\mathbb{E}_P[L_{\mu(x)}(z,X,Y,\theta)|X=x]$ is locally Lipschitz continuous around $\bar{\theta}$. 
\end{proposition}

\noindent
\textbf{Proof.}
Part (i) is established by utilizing Danskin's theorem on $\mathcal{Z}_\beta$ in Proposition \ref{assumption:infcompactness-feasible-set}, see e.g. \cite[Proposition 2.1]{ye2010new} or \cite[Theorem 9.26]{shapiro2021lectures}.
Part (ii) is derived by applying Carath\'eodory's theorem to the convex hull in \eqref{eqn:clarke-subdiff-value-function}.
% $\mathfrak{C}^*(\theta,x) = \min_{z\in \mathcal{Z}}  \mathfrak{C}(z,\theta,x)$.
Part (iii) is obtained directly from part (i).
\hfill $\Box$

In line with 
% Assumption \ref{assumption:inf-compactness} and 
Assumption \ref{assumption:uni-inf-compact}, we adopt the following assumptions.

\begin{assumption} \label{assumption:uni-inf-compact2}
    Let $\theta \in \Theta$ be fixed. For every $x\in \mathcal{X}$, there exists a constant $\delta>0$ such that the set 
    \begin{eqnarray*}
        \Big\{ z: r\in  \nabla_z \mathfrak{C}(z,\theta,x') +{\cal N}_{\cal Z}(z), \mathbb{E}_{P} [L_{\mu(x)}(z,x',Y,\theta)|X=x'] \leq \alpha, r\in \delta\mathcal{B} \Big\}
    \end{eqnarray*}
    is bounded for some constant $\alpha$ and every $x'$ in a close neighborhood of $x$ relative to $\mathcal{X}$.
\end{assumption}

Based on these assumptions, we can establish the necessary optimality condition for problem \eqref{eqn:CP2-mu-vfunc} and moreover problem \eqref{eq:ILO-2-rewrite} by following the scheme in Section \ref{sec:kkt-approach}.

\begin{theorem} \label{thm:first-order-necessary-condition2}
    Suppose that (a) Assumption \ref{assumption:differentiable-regular}-\ref{assumption:lipschitz_in_z}, \ref{assumption:convex-smooth2} %\ref{assumption:stochastic-pc}, %\ref{assumption:infcompactness-feasible-set}, 
    -\ref{assumption:uni-inf-compact2} hold;
    (b) for fixed $\theta$, there exists a nonnegative function $\sigma(x)$ with $\mathbb{E}_P[\sigma(X)]< \infty$ such that 
    \begin{eqnarray*} 
    % \label{eqn:theorem-lipshcitz-bound1}
        % \max\left\{ \| \nabla_\theta\mathbb{E}_P\left[L(z,X,Y)\mid X=x\right] \| , 
        \max\left\{ \nabla_\theta\mathbb{E}_P\left[L_{\mu(x)}(z,X,Y,\theta)\mid X=x\right] \| , \| \nabla_{z\theta}^2 \mathfrak{C}(z,\theta,x) \| \right\}
        % \right\} 
        \leq \sigma(x)
    \end{eqnarray*}
    uniformly for all $z\in \Gamma(\theta,x)$;
    (c) $V_\mu(\theta):=\mathbb{E}_P[v_{\mu(X)}(\theta,X)]$ is well-defined and the Lipschitz modulus of $v_\mu(x)(\theta,x)$ in $\theta$ is bounded by an integrable function $\kappa(x)$.
    If $(\theta,z(\cdot))$ is a local optimal solution of problem \eqref{eq:ILO-2-rewrite}, then
    % Then for a local optimal solution $\theta^*$ of the ILO model \eqref{eq:ILO-3}, 
    there exist a measurable function $\mu(x)$ and a measurable selection of CD multipliers $\eta(x)$ such that 
    \begin{eqnarray}
    \label{eqn:first-order-necessary-condition-convex2}
        \begin{cases}
            0 \in \mathbb{E}_P\left[
            \nabla_\theta \mathbb{E}_P\left[L(z(X),X,Y,\theta)\mid X\right] +
            \mu(X)\left(\nabla_\theta \mathfrak{C}({z}(X),{\theta},X) - \partial_\theta C^*({\theta},X) \right) \right]\\
            \qquad+ \mathbb{E}_P \left[ 
            % \nabla_\theta \mathbb{E}_P\left[L(z(X),X,Y)\mid X\right] + 
            \nabla_{z\theta}^2 \mathfrak{C}(z(X),\theta,X)^{\top} \eta(X) \right] + \mathcal{N}_\Theta(\theta),\\
            0 \in \nabla_z \mathbb{E}_P\left[L(z(x),x,Y,\theta)\mid X=x\right] + \mu(x)\nabla_{z} \mathfrak{C}(z(x),\theta,x) + \nabla_{zz}^2 \mathfrak{C}(z(x),\theta,x)^{\top} \eta(x)  \\
            \qquad+ D^*\N_\mathcal{Z}(z(x), - \nabla_{z}\mathfrak{C}(z(x),\theta,x)) (\eta(x)),  \; \forall \inmat{ a.e. }  x\in \mathcal{X},\\
            0 \in \nabla_z \mathfrak{C}(z(x),\theta,x) +{\cal N}_{\cal Z}(z), \quad  \forall \inmat{ a.e. } x\in \mathcal{X}.
        \end{cases}
    \end{eqnarray}
    % Moreover, if in addition Assumption \ref{assumption:uni-inf-compact2} holds, then theres exist a measurable selection of CD multipliers $\eta(x)$ such that \eqref{eqn:first-order-necessary-condition-convex2} holds.
\end{theorem}

\noindent
\textbf{Proof.} Since $(\theta,z(\cdot))$ is a local optimal solution of problem \eqref{eq:ILO-2-rewrite}, then Assumption \ref{assumption:stochastic-pc} ensures that there exists a measurable function $\mu(x)$ such that $(\theta,z(\cdot))$ is also a local optimal solution of problem \eqref{eqn:CP2-mu-vfunc}.
The rest of proof is similar to Proposition \ref{prop:subdifferentiable-V},
% \cite[Corollary 3.8]{lucet2001sensitivity}, \cite[Proposition 2.15, Theorem 2.16]{xu2010necessary}, 
and 
Theorems \ref{thm:first-order-necessary-condition}-\ref{thm:first-order-necessary-condition-measurable},
we skip the details.
\hfill $\Box$

When $\mathcal{Z}$ is a nonnegative orthant or a polyhedron, the necessary optimality condition can be explicitly written as follows by Proposition \ref{prop:coderivative-first-orthant} and Proposition \ref{prop:normal_to_azb} directly.

\begin{corollary}
    Let the conditions in Theorem \ref{thm:first-order-necessary-condition2} hold. Then we have the following optimality conditions.
    \begin{itemize}
        \item[(i)] If $\mathcal{Z}:= \mathbb{R}_+^{d_z}$, the necessary first-order optimality condition of problem \eqref{eqn:CP} is that there exist $\mu\geq 0$, $\zeta:\mathcal{X}\rightarrow \mathbb{R}^{d_z}$, and $\eta_i(\cdot)\in \mathcal{L}^q(\mathbb{R})$ for $i=1,\cdots,d_z$, 
        such that 
        \begin{eqnarray*}
            \begin{cases}
            0 \in \mathbb{E}_P\left[
            \nabla_\theta \mathbb{E}_P\left[L(z(X),X,Y,\theta)\mid X\right] +
            \mu(X)\left(\nabla_\theta \mathfrak{C}({z}(X),{\theta},X) - \partial_\theta C^*({\theta},X) \right) \right]\\
            \qquad+ \mathbb{E}_P \left[ 
            % \nabla_\theta \mathbb{E}_P\left[L(z(X),X,Y)\mid X\right] + 
            \nabla_{z\theta}^2 \mathfrak{C}(z(X),\theta,X)^{\top} \eta(X) \right] + \mathcal{N}_\Theta(\theta),\\
            0\in \nabla_z \mathbb{E}\left[ L({z}(X),X,Y) |X=x \right] + \mu(x) \nabla_z \mathfrak{C}({z}(x),{\theta},x) + \nabla_{zz}^2 \mathfrak{C} ({z}(x),{\theta},x)^{\top} \eta(x)\\
            \qquad+ \zeta(x),\; \forall \text{ a.e. } x\in \mathcal{X},\\
            \zeta_{L(x)}(x) = 0, \; \eta_{I_+(x)}(x) = 0, \forall \text{ a.e. } x\in \mathcal{X},\\
            \forall i \in I_0(x), \mbox{ either } \zeta_i(x)<0, \eta_i(x)<0 \mbox{ or } \zeta_i(x) \eta_i(x)=0,\; \forall x\in \mathcal{X},\\
            0\in \nabla_z \mathfrak{C}(z(x),\theta,x) +{\cal N}_{\cal Z}(z(x)),\; \forall \text{ a.e. } x\in \mathcal{X},  \\
            \mathfrak{C}(z(x),\theta,x) - \mathfrak{C}^*(\theta,x)\leq 0, \;\forall \text{ a.e. } x\in \mathcal{X},
            \end{cases}
        \end{eqnarray*}
        where 
        \begin{eqnarray*}
L(x)
% :=L(z(x),\nabla_{z}\mathfrak{C}(z(x),\theta,x))
&:=&\{ i\in \{1,2,\dots, d_z\}: z_i(x) >0, (\nabla_{z}\mathfrak{C}(z(x),\theta,x))_i =0\},\\
I_+(x)
% :=I_+(z(x),\nabla_{z}\mathfrak{C}(z(x),\theta,x))
&:=&\{ i\in \{1,2,\dots, d_z\}: z_i(x) =0, (\nabla_{z}\mathfrak{C}(z(x),\theta,x))_i >0\},\\
I_0(x)
% :=I_0(z(x),\nabla_{z}\mathfrak{C}(z(x),\theta,x))
&:=&\{ i\in \{1,2,\dots, d_z\}: z_i(x) =0, (\nabla_{z}\mathfrak{C}(z(x),\theta,x))_i =0\}.
\end{eqnarray*}

    \item[(ii)] If $\mathcal{Z} = \left\{z \in \mathbb{R}_+^{d_z}: Az \leq b\right\}$, the necessary first-order optimality condition of problem \eqref{eqn:CP} is that there exist $\mu\geq 0$, $\zeta:\mathcal{X}\rightarrow \mathbb{R}^{d_z}$, $\eta_i(\cdot)\in \mathcal{L}^q(\mathbb{R})$ for $i=1,\cdots,d_z$,
    $\lambda(\cdot): \mathcal{X} \rightarrow \mathbb{R}_+^{m}$ $J_1(\cdot): \mathcal{X} \rightarrow 2^{\{1,\dots,m\}}$, and $J_2(\cdot): \mathcal{X} \rightarrow 2^{\{1,\dots,m\}}$,
    such that 
\begin{eqnarray*}
            \begin{cases}
            0 \in \mathbb{E}_P\left[
            \nabla_\theta \mathbb{E}_P\left[L(z(X),X,Y,\theta)\mid X\right] +
            \mu(X)\left(\nabla_\theta \mathfrak{C}({z}(X),{\theta},X) - \partial_\theta C^*({\theta},X) \right) \right]\\
            \qquad+ \mathbb{E}_P \left[ 
            % \nabla_\theta \mathbb{E}_P\left[L(z(X),X,Y)\mid X\right] + 
            \nabla_{z\theta}^2 \mathfrak{C}(z(X),\theta,X)^{\top} \eta(X) \right] + \mathcal{N}_\Theta(\theta),\\
            0\in \nabla_z \mathbb{E}\left[ L({z}(X),X,Y) |X=x \right] + \mu(x) \nabla_z \mathfrak{C}({z}(x),{\theta},x) + \nabla_{zz}^2 \mathfrak{C} ({z}(x),{\theta},x)^{\top} \eta(x)\\
            \qquad + \zeta(x), \; \forall \text{ a.e. } x\in \mathcal{X},\\
            - \nabla_{z}\mathfrak{C}(z(x),\theta,x) = A^{\top}\lambda(x),\quad \lambda_i(x)(a_i^{\top}z(x)-b_i) = 0, \quad \forall i=1,\dots,m,\ \text{ a.e. } x\in \mathcal{X},\\
            J_1(x)\subset J_2(x) \subset I_0(z(x)),\; \forall \text{ a.e. } x\in \mathcal{X}, \\
            \zeta(x) = A^{\top}_{J_2(x) \backslash J_1(x)}\rho(x) + A^{\top}_{I_+(z(x))\cup J_1(x)}\nu(x),\; \forall \text{ a.e. } x\in \mathcal{X},\\
            \rho(x)\in \mathbb{R}_+^{|J_2(x)|-|J_1(x)|}, \nu(x)\in \mathbb{R}^{|J_1(x)|+|I_+(z(x))|},\; \forall \text{ a.e. } x\in \mathcal{X},\\
            a_i^{\top}\eta(x) = 0, \forall i\in I_+(z(x))\cup J_1(x),\;  a_i^{\top}\eta(x)\geq 0, \forall i\in J_2(x) \backslash J_1(x),\; \forall \text{ a.e. } x\in \mathcal{X},\\
            0\in \nabla_z \mathfrak{C}(z(x),\theta,x) +{\cal N}_{\cal Z}(z(x)),\; \forall \text{ a.e. } x\in \mathcal{X},  \\
            \mathfrak{C}(z(x),\theta,x) - \mathfrak{C}^*(\theta,x)\leq 0, \;\forall \text{ a.e. } x\in \mathcal{X},
            \end{cases}
        \end{eqnarray*}
    where 
    $I_+(z(x)):= \{ i\in \{1,\dots,m\}: \lambda_i(x)>0, a_i^{\top}z(x) = b_i \}$, and $I_0(z(x)):= \{ i\in \{1,\dots,m\}: \lambda_i(x)=0, a_i^{\top}z(x) = b_i \}$.
    \end{itemize}
    Moreover, for a.e. $x\in \mathcal{X}$, $\partial_\theta C^*({\theta},x)$ can be replaced by $\sum_{i=1}^{d_\theta+1} \lambda_i(x) \nabla_\theta \mathfrak{C}(z_i^*(x),{\theta},x)$ with $z_1^*(x),\dots,z_{d_\theta+1}^*(x)\in \mathcal{Z}^*({\theta},x)$, and $\lambda_1(x),\dots,\lambda_{d_\theta+1}(x)\geq0$ with $\sum_{i=1}^{d_\theta+1} \lambda_i(x)=1$.
\end{corollary}

\noindent
\textbf{Proof.}
By invoking Proposition \ref{prop:coderivative-first-orthant} and Proposition \ref{prop:normal_to_azb} for problem \eqref{eqn:first-order-necessary-condition-convex2}, we derive the corresponding optimality conditions.
Moreover, the finite convex combination representation follows from Proposition~\ref{prop:Cstar-lipschitz-differentiable}~(ii).
\hfill $\Box$

\section{Applications}
\label{sec:motivating-applications}

An ILO problem typically comprises three main components: the decision-making task, a predictive model for the conditional distribution, and the loss function.
In this section, we build upon previously established results to derive the first-order necessary optimality conditions for some specific practical ILO problems.

\subsection{Finite-dimensional parameterized predictive models}

Finite-dimensional parameterized predictive models are prevalent in the literature.
The examples presented below demonstrate the application of finite-dimensional parameterized models within the ILO framework.

\begin{example} \label{example:predictive-model}
We introduce the application of the linear prediction model, the Gaussian mixture model, and the kernel regression model.
\begin{itemize}
    \item[(i)] \textbf{Linear prediction model.} Suppose the cost function is linear in $y$, i.e., $c(z,y) = y^{\top}z$. 
    The conditional objective of the decision-making problem becomes $\mathbb{E}_{P}[c(z,Y)|X=x] = \mathbb{E}_P[Y|X=x]^{\top} z$.
    Consequently, predicting the conditional distribution $f(x)$ effectively reduces to estimating the conditional mean of $Y$ given $X=x$.
    In this setting, Elmachtoub and Grigas \cite{elmachtoub2022smart} consider a linear prediction model parameterized by $\theta\in \mathbb{R}^{d_x\times d_y}$, namely $\mathbb{E}_P[Y|X=x]=\theta^\top x$. Equivalently, the conditional distribution predictor $f_\theta(x)$ can be represented as a Dirac measure concentrated at $\theta^\top x$, i.e., $f_\theta(x)=\delta_{\theta^\top x}$.
    The conditional objective can be written as $\mathbb{E}_{P}[c(z,Y)|X=x] = x^\top \theta z$.

    \item[(ii)] \textbf{Parametric family of densities (Gaussian mixture model).}
    Approximate the true joint distribution $P$ of $(X,Y)$ by a parametric family. 
    A choice proposed by \cite{wang2022robust,yoon2025data} is a $K$-component Gaussian mixture model (GMM) as $\inmat{GM}\left(\left\{p^k,\mu^k,\Sigma^k\right\}_{k=1}^K\right)$ with $p\in \Delta^K$,
    where $p^k$ are mixture weights and $(\mu^k,\Sigma^k)$ are the mean and covariance of the $k$-th Gaussian component.
    The corresponding joint density is 
    \begin{eqnarray*}
        p^{\text{GM}}\Big((x,y);\{p^k,\mu^k,\Sigma^k\}_{k=1}^K\Big) = \sum_{k=1}^K p^k \mathcal{N}((x,y);\mu^k,\Sigma^k)
    \end{eqnarray*}
    with $\mathcal{N}\left(\cdot;\mu^k,\Sigma^k\right)$ denoting the $k$-th multivariate Gaussian density.
    By \cite[Lemma 1]{wang2022robust}, the conditional distribution of $Y$ given $X=x$ is again a Gaussian mixture, $f_\theta(x)=GM\left(\left\{p_{Y|X}^k,\mu_{Y|X}^k,\Sigma_{Y|X}^k\right\}_{k=1}^K\right)$, where,
    by writing the component parameters in block form as 
    \begin{eqnarray*}
        \mu^k := 
\begin{bmatrix} 
\mu_X^k \\ 
\mu_Y^k 
\end{bmatrix} 
\in \mathbb{R}^{d_X+d_Y},\quad
\Sigma^k := 
\begin{bmatrix} 
\Sigma_{XX}^k & \Sigma_{XY}^k \\ 
\Sigma_{YX}^k & \Sigma_{YY}^k 
\end{bmatrix} 
\in \mathbb{S}_+^{d_X+d_Y},
% \quad \text{for} \ k \in [K].
    \end{eqnarray*}
the conditional mean, covariance and mixture weights are
    \begin{eqnarray*}
\mu_{Y|X}^k &=& \mu_Y^k + \Sigma_{YX}^k (\Sigma_{XX}^k)^{-1} (X - \mu_X^k), \\
\Sigma_{Y|X}^k &=& \Sigma_{YY}^k - \Sigma_{YX}^k (\Sigma_{XX}^k)^{-1} \Sigma_{XY}^k, \\
p_{Y|X}^k &=& \frac{p^k \mathcal{N}(X | \mu_X^k, \Sigma_{XX}^k)}{\sum_{j=1}^K p^j \mathcal{N}(X | \mu_X^j, \Sigma_{XX}^j)}.
    \end{eqnarray*}
    Consequently, the predicted objective $\mathbb{E}_{f_\theta(x)}[c(z,Y)]$ can be evaluated by using the above conditional GMM.
    In this case, the DM needs to estimate the GMM parameters $\theta := \{p^k,\mu^k,\Sigma^k\}_{k=1}^K$.

    \item[(iii)] \textbf{Smoothed empirical distribution with kernel method (kernel regression).}
    The Nadaraya–Watson kernel regression is a classical nonparametric approach for estimating conditional quantities from data and has been used in contextual optimization, see e.g. \cite{nadaraya1964estimating,wang2026data,watson1964smooth}.
    Given i.i.d. samples $\{(x_i,y_i)\}_{i=1}^N$ and a kernel $K$ with a bandwidth $\theta$, 
    we form kernel estimators of the joint density and marginal density as 
    % the probability density functions of $(X,Y)$ and $X$ are estimated using kernel density estimation by
    \begin{eqnarray*}
        \hat{p}(x,y;\theta) = \frac{1}{N} \sum_{i=1}^N K_\theta(x-x_i) K_\theta(y-y_i), \quad
        \hat{p}(x;\theta) = \frac{1}{N} \sum_{i=1}^N K_\theta(x-x_i).
    \end{eqnarray*}
    The estimator of the conditional distribution of $Y$ given $X=x$ is
    \begin{eqnarray} \label{eqn:pdf-kernel}
        f_\theta(y;x) =\hat{p}(y|X=x;\theta) = \frac{\hat{p}(x,y;\theta)}{\hat{p}(x;\theta)} = \frac{\sum_{i=1}^N K_\theta(x-x_i) K_\theta(y-y_i)}{\sum_{i=1}^N K_\theta(x-x_i)}.
    \end{eqnarray}
    Accordingly, the conditional objective can be written as 
    \begin{eqnarray*}
        \mathbb{E}[c(z,Y)|X=x;\theta] 
        &=& \frac{ \sum_{i=1}^N K_\theta(x-x_i) \int c(z,Y) K_\theta(y-y_i) dy}{ \sum_{i=1}^N K_\theta(x-x_i)}.
    \end{eqnarray*}
    Under proper conditions, the above reduces to the Nadaraya–Watson form
    \begin{eqnarray*}
        \mathbb{E}[c(z,Y)|X=x;\theta] 
        &=& \frac{ \sum_{i=1}^N K_\theta(x-x_i) c(z,y_i) }{ \sum_{i=1}^N K_\theta(x-x_i)}.
    \end{eqnarray*}
    Finally, the bandwidth $\theta$ needs to be selected to improve the quality of the inducing decisions.
    \end{itemize}
\end{example}

\subsection{Optimality conditions}

Smart predict then optimize (SPO), proposed by Elmachtoub and Grigas \cite{elmachtoub2022smart}, is a prevalent framework for contextual problems in management science and operations research. 
Its core principle is that predictive models should be trained to reduce downstream decision error/regret, rather than focusing solely on statistical prediction error. 
The SPO framework can be naturally incorporated into the integrated learning and optimization (ILO) setting examined in this work.
In \cite{elmachtoub2022smart}, they perform numerical tests of the SPO framework on a portfolio selection problem with a linear prediction model and a newsvendor problem with a kernel regression model.
We also consider these applications and provide the first-order necessary optimality conditions based on the derived theorems.

\subsubsection{The portfolio selection problem with the linear prediction model.} 
\label{section:portfolio}

Consider a portfolio selection problem
\begin{eqnarray*}
    P(r)\qquad\min_{z\in\mathcal{Z}} \mathfrak{C}(z,r):= -r^\top z + \frac{\lambda}{2} z^\top \Sigma z,
\end{eqnarray*}
where $r\in \mathbb{R}^{d_z}$ is the expected return vector of the assets, $\Sigma\in \mathbb{R}^{d_z\times d_z}$ is the positive definite covariance matrix of assets, and the feasible set is $\mathcal{Z}:=\{z\in \mathbb{R}_+^{d_z}:1^\top z\leq 1\}.$
The objective is to find a portfolio that offers the maximum return while minimizing overall variance.
The return vector $r$ is often highly dependent on contextual information $x$, including news, economic factors, social media, and others, while the variances are typically much
more stable and are not as difficult or sensitive to
predict \cite{elmachtoub2022smart}.
Therefore, we mainly focus on the prediction of the return vector in this example.

Denote the optimal value and  the set of optimal solutions to problem $P(r)$ by $\mathfrak{C}^*(r)$ and $\mathcal{Z}^*(r)$, respectively.
Suppose that DM uses a linear prediction model to characterize the relationship between $x$ and $r$, i.e. $\hat{r} = \theta^\top x$ and $\Theta = \mathbb{R}^{d_x\times d_z}$.
Then the optimistic SPO loss for the prediction model in problem $P(r)$ is defined by 
\begin{eqnarray} \label{eqn:SPO-loss}
    L(z,x,r,\theta) &=& \min_{z\in\mathcal{Z}^*(\hat{r})} \mathfrak{C}(z,r) - \mathfrak{C}^*(r)\\
    &=& \min_{z\in\mathcal{Z}^*(\theta^\top x)} \left\{ -r^\top z + \frac{\lambda}{2} z^\top \Sigma z \right\} - \mathfrak{C}^*(r).\nonumber
\end{eqnarray}
This loss function measures the excess cost incurred when making a suboptimal decision due to an imprecise prediction of the parameter $r$.
Note that the first-order optimality condition for problem $P(\hat{r})$ is 
\begin{eqnarray*}
    0\in -\hat{r} + \lambda \Sigma z + \mathcal{N}_{\mathcal{Z}}(z).
\end{eqnarray*}
By substituting each component into problem \eqref{eq:ILO-3} with the corresponding elements defined above and approximating the objective using the SAA method with sample $\{(x_n,r_n)\}_{n=1}^N$, the SPO framework can be reformulated as the following ILO problem:
    \begin{eqnarray*}
        \min_{\theta\in\Theta}
        && 
        % \mathbb{E}_{P_X}\left[\mathbb{E}_{P_{Y|X}}[L(z^*(\theta,x),X,Y)]\right] = 
        \frac{1}{N} \sum_{n=1}^N \left[
        v(\theta,x_n)\right] \\
        \inmat{s.t.} && v(\theta,x_n) :=  \min_{z} \left\{ -r_n^\top z + \frac{\lambda}{2} z^\top \Sigma z \right\} - \mathfrak{C}^*(r_n), \\
        &&\qquad \qquad \quad \inmat{s.t.}\quad 0\in -\theta^\top x_n + \lambda \Sigma z + \mathcal{N}_{\mathcal{Z}}(z), \qquad \forall  n=1,\dots,N.
        % z\in\arg        \min_{z\in \mathcal{Z}} \mathfrak{C}(z,\theta,x), 
        % \; \inmat{for a.e.} . \quad
        % x \in \mathcal{X}.
\end{eqnarray*}
Since $\Sigma$ is positive definite and the linear independence constraint qualification (LICQ) is typically satisfied for $\mathcal{Z}:=\{z\in \mathbb{R}_+^{d_z}:1^\top z\leq 1\}$,  Theorem \ref{thm:first-order-necessary-condition} is applicable.
%can be easily verified.
%Therefore, if 
Let $(\theta^*,\{z_n\}_{n=1}^N)$ is a local optimal solution of the SPO problem. Then 
by Theorem \ref{thm:first-order-necessary-condition}, 
%we have the following first-order necessary optimality condition:
there exist vectors $\eta_1,\dots,\eta_N \in \mathbb{R}^{d_z}$ and $\zeta_1,\dots,\zeta_N \in \mathbb{R}^{d_z}$ such that 
\begin{eqnarray*}
    \begin{cases}
        0= \frac{1}{N} \sum_{n=1}^N x_n \eta_n^\top\\
        0 = -r_n+ \lambda \Sigma z_n + \lambda \Sigma \eta_n + \zeta_n,\quad \forall n=1,\dots,N,\\
        \zeta_n\in D^*\mathcal{N}_\mathcal{Z}(z_n, \theta^\top x_n - \lambda \Sigma z_n)(\eta_n), \quad \forall n=1,\dots,N, \\
        0\in -\theta^\top x_n + \lambda \Sigma z_n + \mathcal{N}_{\mathcal{Z}}(z_n), \quad \forall n=1,\dots,N.
    \end{cases}
\end{eqnarray*}
Moreover,  by Corollary \ref{corollary:simplex}, 
there exist $\beta\in \mathbb{R}^N$, $\eta_1,\dots,\eta_N \in \mathbb{R}^{d_z}$ and $\zeta_1,\dots,\zeta_N  \in \mathbb{R}^{d_z}$, such that 
\begin{eqnarray}\label{eqn:newsvendor-first-order}
    \begin{cases}
        0= \frac{1}{N} \sum_{n=1}^N x_n \eta_n^\top\\
        0 = -r_n+ \lambda \Sigma z_n + \lambda \Sigma \eta_n + \zeta_n,\quad \forall n=1,\dots,N,\\
        \beta_n (1-\boldsymbol{1}^\top z_n) = 0, \; \tau_n(\boldsymbol{1}^\top \eta_n) = 0,\quad \forall n=1,\dots,N, \nonumber\\
        \text{either } \beta_n>0, \boldsymbol{1}^\top \eta_n >0, \text{ or } \beta_n \boldsymbol{1}^\top \eta_n = 0,\quad \forall n=1,\dots,N,\nonumber \\
        (\zeta_n)_{L(n)} = \beta_n, \; (\eta_n)_{I_+(n)}=0,\quad \forall n=1,\dots,N,\nonumber\\
        \forall i\in I_0(n), \text{ either } (\zeta_n)_i<\beta_n, (\eta_n)_i<0, \text{ or } ((\zeta_n)_i - \beta_n)(\eta_n)_i = 0\quad \forall n=1,\dots,N,\nonumber\\
        0\in -\theta^\top x_n + \lambda \Sigma z_n + \mathcal{N}_{\mathcal{Z}}(z_n), \quad \forall n=1,\dots,N,
    \end{cases}
\end{eqnarray}
where $L(n)$, $I_+(n)$, $I_0(n)$, and $\tau_n$ are defined as follows:
for $n=1,\dots,N$,
let $L(n)= \{i=\{1,\dots,d_z\}: (z_n)_i>0\}$;
if $1^\top z_n< 1$, then $I_+(n)=\{i=\{1,\dots,d_z\}: (z_n)_i=0, (\theta^\top x_n - \lambda \Sigma z_n)_i>0\}$, $ I_0(n)=\{i=\{1,\dots,d_z\}: (z_n)_i=0, (\theta^\top x_n - \lambda \Sigma z_n)_i = 0\}$;
otherwise, if $1^\top z_n=1$, then let $\tau_n:=(\theta^\top x_n - \lambda \Sigma z_n)_i$ for any $i\in L(n)$, $I_+(n)=\{i=\{1,\dots,d_z\}: (z_n)_i=0, (\theta^\top x_n - \lambda \Sigma z_n)_i>\tau\}$, $ I_0(n)=\{i=\{1,\dots,d_z\}: (z_n)_i=0, (\theta^\top x_n - \lambda \Sigma z_n)_i = \tau_n\}$.

\subsubsection{The newsvendor problem with a kernel regression.}
\label{section:newsvendor}

Consider a newsvendor problem with a demand distribution predicted as $f_\theta(x)$.  
Let $h>0$ and $b>0$ denote the unit holding cost and the unit backorder cost, respectively.
The total cost is defined as $c(z,y) = h(z-y)_+ + b(y-z)_+, $
where $(a)_+ = \max\{a,0\}$.
The corresponding newsvendor optimization problem is
    \begin{eqnarray} \label{eqn:newsvendor-problem}
    \min_{z\geq 0} \mathbb{E}_{f_\theta(x)}\left[ h(z-Y)_+ + b(Y-z)_+ \right].
    \end{eqnarray}
Denote the optimal value of the problem under the demand distribution $f$ by $\mathfrak{C}^*(f)$.
Observe that if $f$ is a Dirac distribution or a point prediction model, i.e. $f=\delta_{\hat{y}}$, then the optimal solution is simply $z^*(\theta,x)= \hat{y}$.
To obtain a nontrivial solution, we instead use the kernel regression model from Example \ref{example:predictive-model} (iii), which serves as a distributional prediction model for a continuous distribution.
As discussed in Example \ref{example:predictive-model} (iii), the bandwidth $\theta>0$ of the kernel regression model should be selected to improve the quality of the decisions. 
Therefore, we also use the SPO loss \eqref{eqn:SPO-loss} to find the optimal bandwidth.

Denote the cumulative distribution function and the probability density function of $f_\theta(x)$ at $z$ by $F_\theta(z;x)$ and $p_\theta(z;x)$.
Let $\mathfrak{C}(z,\theta,x):=\mathbb{E}_{f_\theta(x)}\left[ h(z-Y)_+ + b(Y-z)_+ \right]$.
The first-order derivative in $z$ of the objective function is $\nabla_z\mathfrak{C}(z,\theta,x)=(h+b)F_\theta(z;x)- b$, and the second-order derivative in $z$ is $\nabla_{zz}^2\mathfrak{C}(z,\theta,x)=(h+b)p_\theta(z;x)$.
The corresponding first-order optimality condition for the newsvendor problem \eqref{eqn:newsvendor-problem} is 
\begin{eqnarray*}
    0\in (h+b) F_\theta(z;x) - b + \mathcal{N}_{\mathbb{R}_+}(z).
\end{eqnarray*}
Given a sample of contextual information $x\in \mathbb{R}^{d_x}$ and historical demand $y\in \mathbb{R}$, $\{(x_n,y_n)\}_{n=1}^N$.
The ILO problem with SPO loss \eqref{eqn:SPO-loss} can be written as
    \begin{eqnarray*}
        \min_{\theta\in\Theta}
        && 
        \frac{1}{N} \sum_{n=1}^N \left[v(\theta,x_n)\right]\\
        \inmat{s.t.} && v(\theta,x_n) :=  \min_{z} \left\{  h(z-y_n)_+ + b(y_n-z)_+ \right\} - \mathfrak{C}^*(y_n),\\
        &&\qquad \qquad \quad \inmat{s.t.}\quad 0\in (h+b) F_\theta(z;x_n) - b + \mathcal{N}_{\mathbb{R}_+}(z), \qquad \forall  n=1,\dots,N.
\end{eqnarray*}
By Corollary \ref{corollary:opt-cond-first-orthant}, we obtain the  first-order necessary optimality condition for the SPO problem, that is, 
if $(\theta^*,\{z_n\}_{n=1}^N)$ is a local optimal solution, then
there exist vectors $\eta=(\eta_1,\dots,\eta_N)^T \in \mathbb{R}^{N}$ and $\zeta=(\zeta_1,\dots,\zeta_N) \in \mathbb{R}^{N}$, such that 
\begin{eqnarray*}
        \begin{cases}
            0 \in  \frac{1}{N} \sum_{n=1}^N (h+b)\nabla_\theta F_\theta(z_n;x_n)\eta_n + \mathcal{N}_\Theta(\theta),\\
            0 \in \partial_z \big(h(z_n- y_n)_+ + b(y_n-z_n)\big) + (h+b)p_\theta(z_n;x_n) \eta_n  + \zeta_n, \forall n =1,\dots,N,\\
            \zeta_{L} = 0, \; \eta_{I_+} = 0, \; \forall n\in I_0, \mbox{ either } \zeta_n<0, \eta_n<0 \mbox{ or } \zeta_n \eta_n=0,\\
            0 \in (h+b) F_\theta(z;x_n) - b +{\cal N}_{\mathbb{R}_+}(z_n), \forall n =1,\dots,N,
        \end{cases}
\end{eqnarray*}
where $L:=\{n\in\{1,\dots,N\}: z_n>0\}$, $I_+:=\{n\in\{1,\dots,N\}: z_n=0, (h+b)F_\theta(z_n;x_n)-b>0 \}$, and $I_0:=\{n\in\{1,\dots,N\}: z_n=0, (h+b)F_\theta(z_n;x_n)-b=0 \}$. 
In particular, $\nabla_z \big(h(z_n- y_n)_+ + b(y_n-z_n)\big)$ is $\{h\}$ if $z_n>y_n$; $\{-b\}$ if $z_n< y_n$; $[-b,h]$ if $z_n=y_n$.

For $\nabla_\theta F_\theta(z;x_n)$, we recall the probability density function given by the kernel regression in \eqref{eqn:pdf-kernel} and the corresponding cumulative distribution function as
\begin{eqnarray*}
    p_\theta(y|x) =\sum_{m=1}^N w_m(x;\theta) K_\theta(y-y_m);\quad
    F_\theta(y;x) 
    % = \frac{\sum_{m=1}^N K_\theta(x-x_m) \int_{-\infty}^y K_\theta(t-y_m)dt}{\sum_{m=1}^N K_\theta(x-x_m)} 
    = \sum_{m=1}^N w_m(x;\theta) \int_{-\infty}^y K_\theta(t-y_m)dt,
\end{eqnarray*}
where $w_m(x;\theta):= \frac{K_\theta(x-x_m)}{\sum_{l=1}^N K_\theta(x-x_l)}$.
When the kernel function is chosen as the Gaussian kernel, i.e., 
\begin{eqnarray*}
    K_\theta(u):= \frac{1}{(2\pi)^{d/2}\theta^d} \exp\left( -\frac{\|u\|^2}{2\theta^2} \right),
\end{eqnarray*}
we have
\begin{eqnarray*}
&&\psi_n(x;\theta):=\nabla_\theta \log K_\theta(x-x_n) = -\frac{d_x}{\theta} + \frac {\|x-x_n\|^2} {\theta^3}, \\ 
    &&\nabla_\theta w_n(x;\theta) = w_n(x;\theta) \left( \psi_n(x;\theta) - \sum_{m=1}^N w_m(x;\theta)\psi_m(x;\theta)\right),\\
    &&\nabla_\theta \int_{-\infty}^y K_\theta(t-y_n)dt = \nabla_\theta \int_{-\infty}^y \frac{1}{ \sqrt{2\pi}\theta} \exp\left( -\frac{(t-y_n)^2}{2\theta^2} \right) dt = \nabla_\theta \Phi\left(\frac{y-y_n}{\theta}\right)\\
    &&\qquad \qquad \qquad\qquad \quad \; \;  = -\phi\left(\frac{y-y_n}{\theta}\right) \frac{y-y_n}{\theta^2} = -\frac{y-y_n}{\theta} K_\theta(y-y_n),
\end{eqnarray*}
where $\Phi$ and $\phi$ is the cumulative distribution function and probability density function of a standard normal distribution, and consequently
\begin{eqnarray*}
    \nabla_\theta F_\theta(y;x) &=& \sum_{m=1}^N  w_m(x;\theta) \left( \psi_m(x;\theta) - \sum_{l=1}^N w_l(x;\theta)\psi_l(x;\theta)\right) \Phi\left(\frac{y-y_m}{\theta}\right)\\
    &&- \sum_{m=1}^N w_m(x;\theta)\frac{y-y_m}{\theta} K_\theta(y-y_m).
\end{eqnarray*}

% \section*{Acknowledgments}
% \bibliographystyle{plain}
\bibliographystyle{apalike}
\bibliography{literature}

\newpage
\begin{appendices}

\section{Appendix}
\label{sec:appendix}

\subsection{Proof of Lemma \ref{lemma:normal_to_gph_normal_polyhedron}}
\label{sec:proof-lemma-3.7}

\noindent
\textbf{Proof of Lemma \ref{lemma:normal_to_gph_normal_polyhedron}.}
To ease the notation, let $G= \text{gph} \mathcal{N}_\mathcal{Z}\subset \mathbb{R}^{d_z} \times \mathbb{R}^{d_z}$.
Since $\mathcal{Z}$ is polyhedral, $\mathcal{N}_\mathcal{Z}$ is a polyhedral set-valued mapping (see \cite[Theorem 4.46]{rockafellar1998variational}) and its graph $G$ can be represented as the union of finitely many polyhedral sets in $\mathbb{R}^{d_z\times d_z}$.
Therefore, by the definition of the limiting normal cone, we have that, for any sufficiently small neighborhood $U$ of $(\bar{z},\bar{v})$, 
\begin{eqnarray} \label{eqn:proof-normal0}
    \mathcal{N}_G(\bar{z},\bar{v}) = \bigcup_{(z,v)\in U\cap G} T_G(z,v)^*.
\end{eqnarray}
For each $(z,v)\in G$, by the Reduction Lemma (see Lemma in \cite{robinson2009local} or Theorem 5.6 in \cite{rockafellar1989proto}), 
the tangent cone to $G$ at $(z,v)$ is 
\begin{eqnarray} \label{eqn:proof-normal1}
    T_G(z,v) = \text{gph} \mathcal{N}_{K_\mathcal{Z}(z,v)}.
\end{eqnarray}
Note that for any convex cone $K$, we have 
\begin{eqnarray*}
    v'\in \mathcal{N}_K (z') \Longleftrightarrow z'\in K, v'\in K^*, z'\perp v', 
\end{eqnarray*}
and thus 
\begin{eqnarray*}
    \text{gph} \mathcal{N}_K = \left\{ (z',v') : z'\in K, v'\in K^*, z'\perp v'  \right\}.
\end{eqnarray*}
Therefore, by \eqref{eqn:proof-normal1},
\begin{eqnarray*}
    T_G(z,v)=\left\{ (z',v') : z'\in K_\mathcal{Z}(z,v),\; v'\in K_\mathcal{Z}(z,v)^*,\; z'\perp v' \right\},
\end{eqnarray*}
and then by the definition of polar cone,
\begin{eqnarray}
    T_G(z,v)^* & = & \left\{ (r,u) : \langle (r,u),(z',v') \rangle \leq 0, \forall (z',v')\in T_G(z,v) \right\}\nonumber\\
    & = & \left\{ (r,u) : \langle r,\z' \rangle + \langle u,v' \rangle \leq 0, \forall z'\in K_\mathcal{Z}(z,v),\; v'\in K_\mathcal{Z}(z,v)^*,\; z'\perp v' \right\}. \label{eqn:proof-normal2}
\end{eqnarray}
In \eqref{eqn:proof-normal2}, we first let $z'=0$, and thus $\langle u,v' \rangle\leq 0$ for all $v'\in K_\mathcal{Z}(z,v)^*$, which means $u\in (K_\mathcal{Z}(z,v)^*)^* = K_\mathcal{Z}(z,v)$.
Similarly, we let $v'=0$, and obtain $r\in K_\mathcal{Z}(z,v)^*$.
Combining these two results, we have 
\begin{eqnarray} \label{eqn:proof-normal3}
    T_G(z,v)^* \subset  K_\mathcal{Z}(z,v)^* \times K_\mathcal{Z}(z,v).
\end{eqnarray}
On the other hand, 
for any fixed $(r,u) \in K_\mathcal{Z}(z,v)^* \times K_\mathcal{Z}(z,v)$,
we have $\langle r,z' \rangle + \langle u,v' \rangle\leq 0 +0\leq 0$ for all $z' \in K_\mathcal{Z}(z,v)$ and $v' \in K_\mathcal{Z}(z,v)^*$, and thus $(r,u)\in T_G(z,v)^*$.
Combining with \eqref{eqn:proof-normal3}, we have 
\begin{eqnarray} \label{eqn:proof-normal4}
    T_G(z,v)^* =  K_\mathcal{Z}(z,v)^* \times K_\mathcal{Z}(z,v).
\end{eqnarray}
Therefore, with \eqref{eqn:proof-normal0},
\begin{eqnarray} \label{eqn:proof-normal3.5}
    \mathcal{N}_G(\bar{z},\bar{v}) = \bigcup_{(z,v)\in U\cap G} K_\mathcal{Z}(z,v)^* \times K_\mathcal{Z}(z,v).
\end{eqnarray}

Next, we characterize all cones $K_\mathcal{Z}(z,v)$ with $(z,v)\in U\cap G$.
For any vector $v\in \mathbb{R}^{d_z}$, let $[v] := \{tv:t\in \mathbb{R}\}$ and $[v]^\perp:=\{ w \in \mathbb{R}^{d_z}: \langle v,w \rangle = 0\}$.
Then for polyhedron $\mathcal{Z}$, 
\begin{eqnarray} \label{eqn:proof-normal4}
    K_\mathcal{Z}(z,v)= T_\mathcal{Z}(z) \cap [v]^\perp.
\end{eqnarray}
Since $\mathcal{Z}$ is a polyhedron, for $z\in\mathcal{Z}$ sufficiently close to $\bar{z}$, we have 
\begin{eqnarray} \label{eqn:proof-normal5}
    T_\mathcal{Z}(z) = T_\mathcal{Z}(\bar{z}) + [z - \bar{z}] \supset T_\mathcal{Z}(\bar{z}),
\end{eqnarray}
and $z-\bar{z}\in T_\mathcal{Z}(\bar{z})$.
For $p\in \mathcal{N}_\mathcal{Z}(z)= T_\mathcal{Z}(z)^*$,
\begin{eqnarray*}
    0 \geq \langle p , q + t(z - \bar{z})\rangle= \langle p , q \rangle + t \langle p,z - \bar{z}\rangle, \quad \forall t\in \mathbb{R}, q\in T_\mathcal{Z}(\bar{z}).  
\end{eqnarray*}
This requires that $\langle p,z - \bar{z}\rangle=0$ (set $t\rightarrow \infty$) and $\langle p, q \rangle \leq 0$ for all $q\in T_\mathcal{Z}(\bar{z})$ (set $t=0$), and thus leads to
\begin{eqnarray} \label{eqn:proof-normal6}
    \mathcal{N}_\mathcal{Z}(z) = \mathcal{N}_\mathcal{Z}(\bar{z}) \cap [z - \bar{z}]^\perp \subset \mathcal{N}_\mathcal{Z}(\bar{z}).
\end{eqnarray}

On the other hand, the closed face of the polyhedron $T_\mathcal{Z}(\bar{z})$ takes the form of $F_\mathcal{Z}(\bar{x};v):= T_\mathcal{Z}(\bar{z})\cap [v]^\perp$ with $v\in \mathcal{N}_\mathcal{Z}(\bar{z})$.
In particular, the critical cone $F_\mathcal{Z}(\bar{x};\bar{v})= K_\mathcal{Z}(\bar{z},\bar{v})$ is itself a closed face of $T_\mathcal{Z}(\bar{z})$.
Consider $v\rightarrow \bar{v}$ with $v\in \mathcal{N}_\mathcal{Z}(\bar{z})$.
Since $[v]^{\perp}\to[\bar{v}]^{\perp}$ and $T_\mathcal{Z}(\bar{z})$ are closed, the outer limit of the faces $F_\mathcal{Z}(\bar{x};v)=T_\mathcal{Z}(\bar{z})\cap[v]^{\perp}$ is contained in
$F_\mathcal{Z}(\bar{x};\bar{v})=T_\mathcal{Z}(\bar{z})\cap[\bar{v}]^{\perp}$. Because $T_\mathcal{Z}(\bar{z})$ has only finitely many closed faces,
this implies that there exists a neighborhood $V$ of $\bar{v}$ in
$\mathcal{N}_\mathcal{Z}(\bar{z})$ such that
\begin{eqnarray} \label{eqn:proof-normal7}
    F_\mathcal{Z}(\bar{x};v)\subset F_\mathcal{Z}(\bar{x};\bar{v}), \forall v\in V.
\end{eqnarray}
Consequently, for all $v$ sufficiently close to $\bar{v}$, $F_\mathcal{Z}(\bar{x};v)$ is not
only a closed face of $T_\mathcal{Z}(\bar{z})$ but also a closed face of $F_\mathcal{Z}(\bar{x};\bar{v})$.

Based on the above discussion, for $(z,v)\in \text{gph}\mathcal{N}_\mathcal{Z}$ sufficiently close to $(\bar{z},\bar{v})$, we have $z = \bar{z}+ z'$ with $z'\in T_\mathcal{Z}(\bar{z})$ and $v\in \mathcal{N}_\mathcal{Z}(z) = \mathcal{N}_\mathcal{Z}(\bar{z}) \cap [z']^\perp \subset \mathcal{N}_\mathcal{Z}(\bar{z})$.
Since $z'\perp v$, we have $[z']\subset [v]^\perp$ and thus
\begin{eqnarray}\label{eqn:proof-normal8}
    K_{\mathcal{Z}}(z,v) &=& (T_\mathcal{Z}(\bar{z}) + [z']) \cap [v]^\perp\nonumber\\
    &=& (T_\mathcal{Z}(\bar{z}) \cap [v]^\perp ) + [z'] = F_\mathcal{Z}(\bar{z};v)+ [z'].
\end{eqnarray}
Let $F_1:= F_\mathcal{Z}(\bar{x};v)$.
Because $v$ is close to $\bar{v}$, we have $F_\mathcal{Z}(\bar{z};v)\subset F_\mathcal{Z}(\bar{z};\bar{v})$ by \eqref{eqn:proof-normal7},
which means $F_\mathcal{Z}(\bar{z};v)$ is a closed face of $F_\mathcal{Z}(\bar{z};\bar{v})$.
Moreover, since $[z']\subset K_\mathcal{Z}(z,v)$ by choosing $0\in F_\mathcal{Z}(\bar{z};v)$ in \eqref{eqn:proof-normal8}, we have $z'\perp v$ and thus $z'\in T_\mathcal{Z}(\bar{z})\cap [v]^\perp = F_\mathcal{Z}(\bar{z};v)$.
We let $F_2$ be a closed face of $F_1$ containing $z'$ in its relative interior, and $F_2$ is also a closed face of $F_\mathcal{Z}(\bar{z};\bar{v})$.
Together with \eqref{eqn:proof-normal8}, we have $K_\mathcal{Z}(z,v) = F_1-F_2$.

Conversely, for any given $F_2\subset F_1 \subset F_\mathcal{Z}(\bar{z};\bar{v})$, there exists $v\in \mathcal{N}_\mathcal{Z}(\bar{z})$ such that $F_1 = T_\mathcal{Z}(\bar{z})\cap [v]^\perp$.
Then we choose $z'\in F_2$ sufficiently small such that $z := \bar{z}+z'\in \mathcal{Z}$.
Since $z'\in F_1 = T_\mathcal{Z}(\bar{z})\cap [v]^\perp$, we have $z'\perp v$, and thus $v\in \mathcal{N}_\mathcal{Z}(\bar{z})\cap [z']^\perp$.
Using \eqref{eqn:proof-normal5} and \eqref{eqn:proof-normal6} again, we have 
\begin{eqnarray*}
    T_\mathcal{Z}(z) = T_\mathcal{Z}(\bar{z}) + [z'], \qquad \mathcal{N}_\mathcal{Z}(z) = \mathcal{N}_\mathcal{Z}(\bar{z})\cap  [z']^\perp,
\end{eqnarray*}
and thus $v\in \mathcal{N}_\mathcal{Z}(z)$.
Therefore, 
$(z,v)\in \text{gph}\mathcal{N}_\mathcal{Z}$ and 
\begin{eqnarray*}
    F_1 - F_2 &=& F_1+ [z'] = (T_\mathcal{Z}(\bar{z})\cap [v]^\perp) + [z'] \\
    &=& (T_\mathcal{Z}(\bar{z}) + [z']) \cap [v]^\perp= T_\mathcal{Z}(z) \cap [v]^\perp = K_\mathcal{Z}(z,v).
\end{eqnarray*}

Now we have proved that $K_\mathcal{Z}(z,v)$ arising at points of $(z,v)\in U\cap G$ are exactly the cones $F_1-F_2$ with $F_2\subset F_1$ being closed faces of $K_\mathcal{Z}(\bar{z},\bar{v})$. 
Take this form back into \eqref{eqn:proof-normal3.5}, we arrive at the result.
% \eqref{eqn:normal_to_gph_normal_polyhedron}.
\hfill $\Box$

\subsection{Proof of Corollary \ref{corollary:simplex}}
\label{sec:proof-corollary-5.1}

\noindent
\textbf{Proof of Corollary \ref{corollary:simplex}.}
The feasible set can be written in the form $Az\leq b$ with $m=d_z+1$ in Proposition \ref{prop:normal_to_azb} by setting
\begin{eqnarray*}
    A= 
    \begin{bmatrix}
        -I^{d_z\times d_z}\\
        \boldsymbol{1}^\top
    \end{bmatrix}, \qquad
    b=
    \begin{bmatrix}
        \boldsymbol{0}\\
        1
    \end{bmatrix},
\end{eqnarray*}
where $\boldsymbol{0},\boldsymbol{1}$ are $d_z$-dimensional vectors.
Then for multiplier $(\lambda_1,\dots,\lambda_{d_z},\tau)^\top\in \mathbb{R}_+^{d_z+1}$, 
\begin{eqnarray*}
    -v=A^\top \lambda = -\lambda_{1:d_z} + \tau \boldsymbol{1},\quad \tau(\boldsymbol{1}^\top z-1)=0,\quad \lambda_i z_i = 0,\ \forall i=1,\dots,d_z.
\end{eqnarray*}
Let $L:= \{i\in \{1,\dots,d_z\}: z_i>0,\lambda_i=0\}$, $I_+=\{i\in \{1,\dots,d_z\}: z_i=0, \lambda_i>0\}$, and $I_0=\{i\in \{1,\dots,d_z\}: z_i=0, \lambda_i=0\}$.

Case 1: $\boldsymbol{1}^\top z-1 <0$ and $\tau=0$. The constraints reduce to the nonnegative orthant discussed in Corollary \ref{prop:coderivative-first-orthant} and $\mathcal{N}_{\text{gph}\mathcal{N}_\mathcal{Z}}(z,-v)$ is same as \eqref{eqn:coderivative-first-orthant}.

Case 2: $\boldsymbol{1}^\top z-1 =0$ and $\tau>0$. Then by Proposition \ref{prop:normal_to_azb},
\begin{eqnarray*}
        \mathcal{N}_{\text{gph}\mathcal{N}_\mathcal{Z}}(z,-v)&=&\bigg\{(\zeta,-\eta)\in \mathbb{R}^{d_z}\times \mathbb{R}^{d_z}: \exists (\lambda_1,\dots,\lambda_{d_z},\tau)^\top\in \mathbb{R}_+^{d_z+1} \text{ and } J_1\subset J_2 \subset \{1,\dots,d_z\}, \nonumber\\
    & &\quad -v_i  = -\lambda_i +\tau,\quad \lambda_i z_i= 0, \quad \forall i=1,\dots,d_z,\nonumber\\
    & &\quad \boldsymbol{1}^\top z-1 = 0, \ \tau>0,  \nonumber\\
    & &\quad z_i = 0,\; \lambda_i=0,\; \forall i\in J_2,\nonumber\\
    & &\quad \zeta = \big(-I^{d_z\times d_z}\big)^{\top}_{J_2 \backslash J_1}\mu + \big(-I^{d_z\times d_z}\big)^{\top}_{\{i: z_i = 0,\lambda_i>0\}\cup J_1}\nu + \boldsymbol{1}^\top \beta,\\
    && \quad\mu\in \mathbb{R}_+^{|J_2|-|J_1|}, \nu\in \mathbb{R}^{|J_1|+|\{i: z_i = 0, \lambda_i>0\}|}, \beta\in \mathbb{R},\nonumber\\
    & &\quad \eta_i = 0, \forall i\in \{i: z_i = 0, \lambda_i>0\}\cup J_1 \cup \{d_z+1\},\;  \eta_i\leq 0, \forall i\in J_2 \backslash J_1 \bigg\}.
    \end{eqnarray*}
Note that $\lambda_i=v_i+\tau$, we have $L= \{i: z_i>0,v_i =- \tau\}$, $I_+=\{i: z_i=0, v_i>-\tau\}$, and $I_0=\{i: z_i=0, v_i = -\tau\}$. 
Following the same discussion in the proof of Corollary \ref{prop:coderivative-first-orthant}, then 
\begin{eqnarray*}
        \mathcal{N}_{\text{gph}\mathcal{N}_\mathcal{Z}}(z,-v)&=&\bigg\{(\zeta,-\eta)\in \mathbb{R}^{d_z}\times \mathbb{R}^{d_z}: \exists \beta\in \mathbb{R},\\
        &&\quad \zeta_L = \beta, \; \eta_{I_+}=0, \; \boldsymbol{1}^\top \eta=0,\\
        &&\quad \forall i\in I_0, \text{ either } \zeta_i<\beta, \eta_i<0, \text{ or } (\zeta_i - \beta)\eta_i = 0 \bigg\}.
    \end{eqnarray*}

Case 3: $\boldsymbol{1}^\top z-1 =0$ and $\tau=0$. Then we have
\begin{eqnarray*}
        \mathcal{N}_{\text{gph}\mathcal{N}_\mathcal{Z}}(z,-v)&=&\bigg\{(\zeta,-\eta)\in \mathbb{R}^{d_z}\times \mathbb{R}^{d_z}: \exists (\lambda_1,\dots,\lambda_{d_z},\tau)^\top\in \mathbb{R}_+^{d_z+1} \text{ and } J_1\subset J_2 \subset \{1,\dots,d_z,d_z+1\}, \nonumber\\
    & &\quad -v_i  = -\lambda_i,\quad \lambda_i z_i= 0, \quad \forall i=1,\dots,d_z,\nonumber\\
    & &\quad \boldsymbol{1}^\top z-1 = 0, \ \tau=0,  \nonumber\\
    & &\quad z_i = 0,\; \lambda_i=0,\; \forall i\in J_2\backslash \{d_z+1\},\nonumber\\
    & &\quad \zeta = A^{\top}_{J_2 \backslash J_1}\mu + A^{\top}_{I_+\cup J_1}\nu,\quad\mu\in \mathbb{R}_+^{|J_2|-|J_1|}, \nu\in \mathbb{R}^{|J_1|+|I_+|},\nonumber\\
    & &\quad a_i^\top \eta = 0, \forall i\in I_+\cup J_1,\;  a_i^\top\eta\geq 0, \forall i\in J_2 \backslash J_1 \bigg\}.
    \end{eqnarray*}
% \section{Sample average approximation}
% {\color{red}
% Convergence of stationary points of the SAA problem}
Therefore, $L:= \{i\in \{1,\dots,d_z\}: z_i>0,v_i=0\}$, $I_+=\{i\in \{1,\dots,d_z\}: z_i=0, v_i>0\}$, and $I_0=\{i\in \{1,\dots,d_z\}: z_i=0, v_i=0\}$.
Let 
\begin{eqnarray*}
    \beta:= 
    \begin{cases}
        0, & d_z+1\notin J_2,\\
        v_{d_z+1}\in \mathbb{R}, & d_z+1\in J_1,\\
        \mu_{d_z+1}\in \mathbb{R}, & d_z+1\in J_2\backslash J_1.
    \end{cases}
\end{eqnarray*}
If $d_z+1\in J_1$, then $\boldsymbol{1}^\top \eta=0$; and if $d_z+1 \in J_2\backslash J_1$, then $\boldsymbol{1}^\top \eta\geq 0$, which means either $\beta>0$, $\boldsymbol{1}^\top \eta>0$ or $\beta \boldsymbol{1}^T\eta=0$.
On the other hand, as in Case 2, we have $\zeta_L=\beta$;
$\zeta_{I_+}\in \mathbb{R}$ is free and $\eta_{I_+}=0$.
If $i\in J_1\cap \{1,\dots,d_z\}$, then $\eta_i=0$; and
if $i\in (J_2\backslash J_1)\cap \{1,\dots,d_z\}$, then $\eta_i\leq 0$.
Therefore, for $i\in I_0$, either $\zeta_i < \beta$, $\eta_i <0$ or $(\zeta_i- \beta)\eta_i = 0$.
Then
\begin{eqnarray*}
        \mathcal{N}_{\text{gph}\mathcal{N}_\mathcal{Z}}(z,-v)&=&\bigg\{(\zeta,-\eta)\in \mathbb{R}^{d_z}\times \mathbb{R}^{d_z}: \exists \beta\in \mathbb{R},\nonumber\\
        && \quad \text{ either } \beta>0, \boldsymbol{1}^\top \eta >0, \text{ or } \beta \boldsymbol{1}^\top \eta = 0, \nonumber\\
        &&\quad \zeta_L = \beta, \; \eta_{I_+}=0,\nonumber\\
        &&\quad \forall i\in I_0, \text{ either } \zeta_i<\beta, \eta_i<0, \text{ or } (\zeta_i - \beta)\eta_i = 0 \bigg\}.
    \end{eqnarray*}
Consequently, we summarize the three cases and arrive at the result.
\hfill $\Box$

\end{appendices}

% \section{An example appendix} 

% \section*{Appendix}

\end{document}